\documentclass[12pt]{article}
\usepackage{amssymb}
\newtheorem{theorem}{Theorem}[section]
\newtheorem{proposition}[theorem]{Proposition}
\newtheorem{lemma}[theorem]{Lemma}
\newtheorem{corollary}[theorem]{Corollary}
\def\bull{\vrule height .9ex width .8ex depth -.1ex}
\newenvironment{proof}{\smallbreak \noindent {\bf Proof.~}}
              {\unskip\nobreak\hfill\hskip 2em \bull\par\medbreak}
\newenvironment{proofof}[1]{\medbreak\noindent{\bf Proof of~#1.~}}
              {\unskip\nobreak\hfill\hskip 2em \bull\par\medbreak}
\def\bZ{\mathbb{Z}}
\def\cG{\mathcal{G}}
\def\tQ{\widetilde Q}
\def\tG{\widetilde G}
\def\hphi{\widehat\phi}
\def\hpsi{\widehat\psi}
\def\hD{\widehat D}
\def\eps{\varepsilon}
\def\si{\sigma}
\def\la{\lambda}
\def\Aut{\mathop{\mathrm{Aut}}}

\newcounter{thesame}
\title{On a series of finite automata\\defining free transformation groups}
\author{Mariya Vorobets\thanks{%
        Partially supported by the NSF grants DMS-0308985 and DMS-0456185.}
        {} and Yaroslav Vorobets$^\fnsymbol{thesame}$\thanks{%
        The second author is supported by a Clay Research Scholarship.}
}
\date{}
\begin{document}

\maketitle

\begin{abstract}
We introduce two series of finite automata starting from the so-called
Aleshin and Bellaterra automata.  We prove that each automaton in the first
series defines a free non-Abelian group while each automaton in the second
series defines the free product of groups of order $2$.  Furthermore, these
properties are shared by disjoint unions of any number of distinct automata
from either series.
\end{abstract}

\section{Introduction}\label{main}

A (Mealy) automaton over a finite alphabet $X$ is determined by the set of
internal states, the state transition function and the output function.  A
finite (or finite-state) automaton has finitely many internal states.  An
initial automaton has a distinguished initial state.  Any initial automaton
over $X$ defines a transformation $T$ of the set $X^*$ of finite words in
the alphabet $X$.  That is, the automaton transduces any input word $w\in
X^*$ into the output word $T(w)$.  The transformation $T$ preserves the
lengths of words and common beginnings.  The set $X^*$ is endowed with the
structure of a regular rooted tree so that $T$ is an endomorphism of the
tree.  A detailed account of the theory of Mealy automata is given in
\cite{GNS}.

The set of all endomorphisms of the regular rooted tree $X^*$ is of
continuum cardinality.  Any endomorphism can be defined by an automaton.
However the most interesting are finite automaton transformations that
constitute a countable subset.  If $T_1$ and $T_2$ are mappings defined by
finite initial automata over the same alphabet $X$, then their composition
is also defined by a finite automaton over $X$.  If a finite automaton
transformation $T$ is invertible, then the inverse transformation is also
defined by a finite automaton.  Furthermore, there are simple algorithms to
construct the corresponding composition automaton and inverse automaton.
In particular, all invertible transformations defined by finite automata
over $X$ constitute a transformation group $\cG(X)$.  This fact was
probably first observed by Ho\v{r}ej\v{s} \cite{H}.

A finite non-initial automaton $A$ over an alphabet $X$ defines a finite
collection of transformations of $X^*$ corresponding to various choices of
the initial state.  Assuming all of them are invertible, these
transformations generate a group $G(A)$, which is a finitely generated
subgroup of $\cG(X)$.  We say that the group $G(A)$ is defined by the
automaton $A$.  The groups defined by finite automata were introduced by
Grigorchuk \cite{G} in connection with the Grigorchuk group of intermediate
growth.  The finite automaton nature of this group has great impact on
its properties.  The formalization of these properties has resulted in the
notions of a branch group (see \cite{BGS}), a fractal group (see
\cite{BGN}), and, finally, the most general notion of a self-similar group
\cite{N}, which covers all automaton groups.

The main issue of this paper are free non-Abelian groups of finite
automaton transformations.  Also, we are interested in the free products of
groups of order $2$ (such a product contains a free subgroup of index $2$).
Brunner and Sidki \cite{BS} proved that the free group embeds into the
group of finite automaton transformations over a $4$-letter alphabet.
Olijnyk \cite{O1}, \cite{O2} showed that the group of finite automaton
transformations over a $2$-letter alphabet contains a free group as well as
free products of groups of order $2$.  In the above examples, all automata
are of linear algebraic origin.

A harder problem is to present the free group as the group defined by a
single finite non-initial automaton.  This problem was solved by Glasner
and Mozes \cite{GM}.  They constructed infinitely many finite automata of
algebraic origin that define transformation groups with various
properties, in particular, free groups.  A finite automaton that defines
the free product of $3$ groups of order $2$ was found by Muntyan and
Savchuk (see \cite{N} and Theorem \ref{main4} below).

\begin{figure}[t]
%
\includegraphics{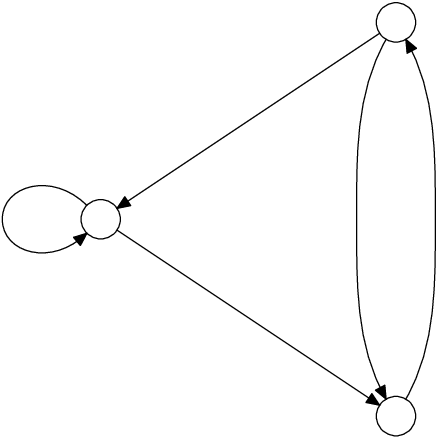}
\includegraphics{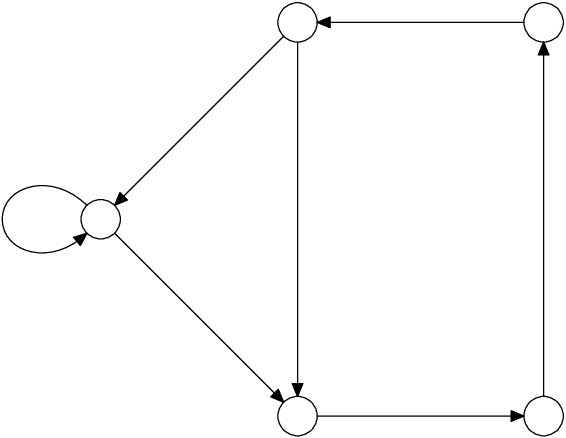}
\begin{picture}(0,0)(0,85)
\put(20,0)
{
 \begin{picture}(0,0)(0,0)
 \put(119,66){$a$}
 \put(34.5,-17){$b$}
 \put(119,-71){$c$}
 \put(66,33){\footnotesize $1|0$}
 \put(66,-39){\footnotesize $1|0$}
 \put(14,15){\footnotesize $0|1$}
 \put(95,-3){\footnotesize $0|1$}
 \put(131.5,-3){\footnotesize $\begin{array}{c} 0|0 \\ 1|1 \end{array}$}
 \end{picture}
}
\put(195.5,0)
{
 \begin{picture}(0,0)(0,0)
 \put(90.5,66){$a$}
 \put(34.5,-17){$b$}
 \put(90.5,-71){$c$}
 \put(161,-73){$d$}
 \put(162,66){$e$}
 \put(53,33){\footnotesize $1|0$}
 \put(53,-39){\footnotesize $1|0$}
 \put(14,15){\footnotesize $0|1$}
 \put(77.5,-3){\footnotesize $0|1$}
 \put(143,-3){\footnotesize $\begin{array}{c} 0|0 \\ 1|1 \end{array}$}
 \put(118,-45.5){\footnotesize $\begin{array}{c} 0|0 \\ 1|1 \end{array}$}
 \put(118,39.5){\footnotesize $\begin{array}{c} 0|0 \\ 1|1 \end{array}$}
 \end{picture}
}
\end{picture}
\vspace{160bp}
\caption{
\label{fig1}
Aleshin's automata.
}
\end{figure}

Actually, the first attempt to embed the free non-Abelian group into a
group of finite automaton transformations was made by Aleshin \cite{A} a
long ago.  He introduced two finite initial automata over alphabet
$\{0,1\}$ and claimed that two automorphisms of the rooted binary tree
$\{0,1\}^*$ defined by these automata generate a free group.  However the
argument in \cite{A} seems to be incomplete.  Aleshin's automata are
depicted in Figure \ref{fig1} by means of Moore diagrams.  The Moore
diagram of an automaton is a directed graph with labeled edges.  The
vertices are the states of the automaton and edges are state transition
routes.  Each label consists of two letters from the alphabet.  The left
one is the input field, it is used to choose a transition route.  The right
one is the output generated by the automaton.  Aleshin considered these
automata as initial, with initial state $b$.

The Aleshin automata are examples of bi-reversible automata.  This notion,
which generalizes the notion of invertibility, was introduced in
\cite{MNS} (see also \cite{GM}).  The class of bi-reversible automata is in
a sense opposite to the class of automata defining branch groups.  All
automata considered in this paper are bi-reversible.

In this paper, we are looking for finite automata that define free
non-Abelian groups of maximal rank, i.e., the free rank of the group is
equal to the number of states of the automaton.  Note that the automata
constructed by Glasner and Mozes do not enjoy this property.  For any of
those automata, the transformations assigned to various internal states
form a symmetric generating set so that the free rank of the group is half
of the number of the states.  Brunner and Sidki conjectured (see \cite{S})
that the first of two Aleshin's automata shown in Figure \ref{fig1} is the
required one.  The conjecture was proved in \cite{VV}.

\begin{theorem}[\cite{VV}]\label{main1}
The first Aleshin automaton defines a free group on $3$ generators.
\end{theorem}

In this paper we generalize and extend Theorem \ref{main1} in several
directions.

The two automata of Aleshin are related as follows.  When the first
automaton is in the state $c$, it is going to make transition to the state
$a$ independently of the next input letter, which is sent directly to the
output.  The second automaton is obtained from the first one by inserting
two additional states on the route from $c$ to $a$ (see Figure \ref{fig1}).

For any integer $n\ge1$ we define a $(2n+1)$-state automaton $A^{(n)}$ of
Aleshin type.  Up to renaming of internal states, $A^{(n)}$ is obtained
from the first Aleshin automaton by inserting $2n-2$ additional states on
the route from $c$ to $a$ (for a precise definition, see Section
\ref{series}); in particular, $A^{(1)}$ and $A^{(2)}$ are the Aleshin
automata.  The Moore diagram of the automaton $A^{(3)}$ is depicted in
Figure \ref{fig6} below.  Note that the number of internal states of an
Aleshin type automaton is always odd.  This is crucial for the proof of the
following theorem.

\begin{theorem}\label{main2}
For any $n\ge1$ the automaton $A^{(n)}$ defines a free group on $2n+1$
generators.
\end{theorem}

Given a finite number of automata $Y^{(1)},\dots,Y^{(k)}$ over the same
alphabet with disjoint sets of internal states $S_1,\dots,S_k$, we can
regard them as a single automaton $Y$ with the set of internal states
$S_1\cup\dots\cup S_k$.  The automaton $Y$ is called the disjoint union of
the automata $Y^{(1)},\dots,Y^{(k)}$ as its Moore diagram is the disjoint
union of the Moore diagrams of $Y^{(1)},\dots,Y^{(k)}$.  The group defined
by $Y$ is generated by the groups $G(Y^{(1)}),\dots,G(Y^{(k)})$.

We define the Aleshin type automata so that their sets of internal states
are disjoint.  Hence the disjoint union of any finite number of distinct
automata of Aleshin type is well defined.

\begin{theorem}\label{main3}
Let $N$ be a nonempty set of positive integers and denote by $A^{(N)}$
the disjoint union of automata $A^{(n)}$, $n\in N$.  Then the automaton
$A^{(N)}$ defines a free group on $\sum_{n\in N}(2n+1)$ generators.
\end{theorem}

One consequence of Theorem \ref{main3} is that the $8$ transformations
defined by the two Aleshin automata generate a free group on $8$
generators.  In particular, any two of them generate a free non-Abelian
group.  Thus Aleshin's claim is finally justified.

\begin{figure}[t]
%
\includegraphics{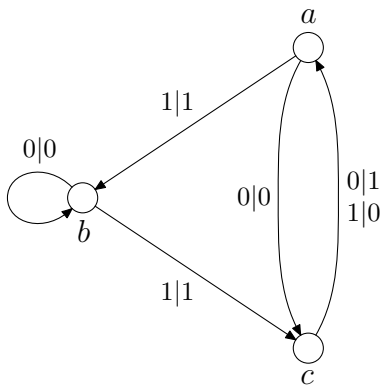}
\begin{picture}(0,0)(-110,85)
 \put(119,66){$a$}
 \put(34.5,-17){$b$}
 \put(119,-71){$c$}
 \put(66,33){\footnotesize $1|1$}
 \put(66,-39){\footnotesize $1|1$}
 \put(14,15){\footnotesize $0|0$}
 \put(95,-3){\footnotesize $0|0$}
 \put(131.5,-3){\footnotesize $\begin{array}{c} 0|1 \\ 1|0 \end{array}$}
\end{picture}
\vspace{160bp}
\caption{
\label{fig2}
The Bellaterra automaton.
}
\end{figure}

The Bellaterra automaton $B$ is a $3$-state automaton over a $2$-letter
alphabet.  Its Moore diagram is depicted in Figure \ref{fig2}.  The
automaton $B$ coincides with its inverse automaton and hence all $3$
transformations defined by $B$ are involutions.  Otherwise there are no
more relations in the group $G(B)$.

\begin{theorem}[\cite{N}]\label{main4}
The Bellaterra automaton defines the free product of $3$ groups of order
$2$.
\end{theorem}

Theorem \ref{main4} is due to Muntyan and Savchuk.  It was proved during
the 2004 summer school on automata groups at the Autonomous University of
Barcelona and so the automaton $B$ was named after the location of the
university.

The Bellaterra automaton $B$ is closely related to the Aleshin automaton
$A$.  Namely, the two automata share the alphabet, internal states, and
the state transition function while their output functions never coincide.
We use this relation to define a series $B^{(1)},B^{(2)},\dots$ of automata
of Bellaterra type.  By definition, $B^{(n)}$ is a $(2n+1)$-state automaton
obtained from $A^{(n)}$ by changing values of the output function at all
elements of its domain.  Also, we define a one-state automaton $B^{(0)}$
that interchanges letters $0$ and $1$ of the alphabet.  All transformations
defined by a Bellaterra type automaton are involutions.

\begin{theorem}\label{main5}
For any $n\ge0$ the automaton $B^{(n)}$ defines the free product of $2n+1$
groups of order $2$.
\end{theorem}

\begin{theorem}\label{main6}
Let $N$ be a nonempty set of nonnegative integers and denote by $B^{(N)}$
the disjoint union of automata $B^{(n)}$, $n\in N$.  Then the automaton
$B^{(N)}$ defines the free product of $\sum_{n\in N}(2n+1)$ groups of
order $2$.
\end{theorem}

Theorems \ref{main3} and \ref{main6} have the following obvious corollary.

\begin{corollary}\label{main7}
(i) Let $n$ be an integer such that $n=3$ or $n=5$ or $n\ge7$.  Then there
exists an $n$-state automaton over alphabet $\{0,1\}$ that define a free
transformation group on $n$ generators.

(ii) For any integer $n\ge3$ there exists an $n$-state automaton over
alphabet $\{0,1\}$ that define a transformation group freely generated by
$n$ involutions.
\end{corollary}

We prove Theorems \ref{main1}, \ref{main2}, and \ref{main3} using the dual
automaton approach.  Namely, each finite automaton $Y$ is assigned a dual
automaton $Y'$ obtained from $Y$ by interchanging the alphabet with the set
of internal states and the state transition function with the output
function.  It turns out that there is a connection between transformation
groups defined by $Y$ and $Y'$.  As intermediate results, we obtain some
information on the dual automata of the Aleshin type automata.

\begin{proposition}\label{main8}
(i) The dual automaton of the Aleshin automaton defines a group that acts
transitively on each level of the rooted ternary tree $\{a,b,c\}^*$.

(ii) For any $n\ge1$ the dual automaton of $A^{(n)}$ defines a group that
acts transitively on each level of the rooted $(2n+1)$-regular tree
$Q_n^*$.
\end{proposition}

The proof of Theorem \ref{main4} given in \cite{N} also relies on the dual
automaton approach.  In particular, it involves a statement on the dual
automaton $\hD$ of $B$.  Since the group $G(B)$ is generated by
involutions, it follows that the set of double letter words over the
alphabet $\{a,b,c\}$ is invariant under the action of the group $G(\hD)$.
Hence $G(\hD)$ does not act transitively on levels of the rooted tree
$\{a,b,c\}$.

\begin{proposition}[\cite{N}]\label{main9}
The dual automaton of the Bellaterra automaton defines a transformation
group that acts transitively on each level of the rooted subtree of
$\{a,b,c\}^*$ formed by no-double-letter words.
\end{proposition}

We derive Theorems \ref{main4}, \ref{main5}, and \ref{main6} from Theorem
\ref{main3}.  This does not involve dual automata.  Nonetheless we obtain a
new proof of Proposition \ref{main9} that also works for all Bellaterra
type automata.

\begin{proposition}\label{main10}
For any $n\ge1$ the dual automaton of $B^{(n)}$ defines a group that
acts transitively on each level of the rooted subtree of $Q_n^*$ formed
by no-double-letter words.
\end{proposition}

Finally, we establish relations between groups defined by automata of
Aleshin type and of Bellaterra type.

\begin{proposition}\label{main11}
(i) The group $G(A)$ is an index $2$ subgroup of $G(B^{(\{0,1\})})$;

(ii) for any $n\ge1$ the group $G(A^{(n)})$ is an index $2$ subgroup of
$G(B^{(\{0,n\})})$;

(iii) for any nonempty set $N$ of positive integers the group $G(A^{(N)})$
is an index $2$ subgroup of $G(B^{(N\cup\{0\})})$.
\end{proposition}

\begin{proposition}\label{main12}
(i) $G(A)\cap G(B)$ is a free group on $2$ generators and an index $2$
subgroup of $G(B)$.

(ii) For any $n\ge1$, $G(A^{(n)})\cap G(B^{(n)})$ is a free group on $2n$
generators and an index $2$ subgroup of $G(B^{(n)})$.

(ii) For any nonempty set $N$ of positive integers, $G(A^{(N)})\cap
G(B^{(N)})$ is an index $2$ subgroup of $G(B^{(N)})$.  Also,
$G(A^{(N)})\cap G(B^{(N)})$ is a free group of rank less by $1$ than the
free rank of $G(A^{(N)})$.
\end{proposition}

The paper is organized as follows.  Section \ref{auto} addresses some
general constructions concerning automata and their properties.  In Section
\ref{a} we recall constructions and arguments of the paper \cite{VV} where
Theorem \ref{main1} was proved.  In Section \ref{series} they are applied
to the Aleshin type automata, which results in the proof of Theorem
\ref{main2} (Theorem \ref{series7}).  Besides, Proposition \ref{main8} is
established in Sections \ref{a} and \ref{series} (see Corollaries
\ref{a5plus} and \ref{series6plus}).  In Section \ref{union} we consider
disjoint unions of Aleshin type automata and obtain Theorem \ref{main3}
(Theorem \ref{union6}).  Section \ref{b} is devoted to the study of the
Bellaterra automaton, automata of Bellaterra type, and their relation to
automata of Aleshin type.  Here we prove Theorems \ref{main4}, \ref{main5},
and \ref{main6} (Theorems \ref{b3} and \ref{b4}), Propositions \ref{main9}
and \ref{main10} (Propositions \ref{b9} and \ref{b10}), Proposition
\ref{main11} (Proposition \ref{b2}), and Proposition \ref{main12}
(Propositions \ref{b6}, \ref{b7}, and \ref{b8}).

\section{Automata}\label{auto}

An {\em automaton\/} $A$ is a quadruple $(Q,X,\phi,\psi)$ formed by two
nonempty sets $Q$ and $X$ along with two maps $\phi:Q\times X\to Q$ and
$\psi:Q\times X\to X$.  The set $X$ is to be finite, it is called the {\em
(input/output) alphabet\/} of the automaton.  We say that $A$ is an
automaton over the alphabet $X$.  $Q$ is called the set of {\em internal
states\/} of $A$.  The automaton $A$ is called {\em finite\/} (or {\em
finite-state\/}) if the set $Q$ is finite.  $\phi$ and $\psi$ are called
the {\em state transition function\/} and the {\em output function},
respectively.  One may regard these functions as a single map $(\phi,\psi):
Q\times X\to Q\times X$.

The automaton $A$ canonically defines a collection of transformations.
First we introduce the set on which these transformations act.  This is the
set of words over the alphabet $X$, which is denoted by $X^*$.  A {\em
word\/} $w\in X^*$ is merely a finite sequence whose elements belong to
$X$.  The elements of $w$ are called {\em letters\/} and $w$ is usually
written so that its elements are not separated by delimiters.  The number
of letters of $w$ is called its {\em length\/}.  It is assumed that $X^*$
contains the empty word $\varnothing$.  The set $X$ is embedded in $X^*$ as
the subset of one-letter words.  If $w_1=x_1\dots x_n$ and $w_2=y_1\dots
y_m$ are words over the alphabet $X$ then $w_1w_2$ denotes their
concatenation $x_1\dots x_ny_1\dots y_m$.  The operation $(w_1,w_2)\mapsto
w_1w_2$ makes $X^*$ into the free monoid generated by all elements of $X$.
The unit element of the monoid $X^*$ is the empty word.  Another structure
on $X^*$ is that of a rooted $k$-regular tree, where $k$ is the cardinality
of $X$.  Namely, we consider a graph with the set of vertices $X^*$ where
two vertices $w_1,w_2\in X^*$ are joined by an edge if $w_1=w_2x$ or
$w_2=w_1x$ for some $x\in X$.  The root of the tree is the empty word.  For
any integer $n\ge0$ the $n$-th {\em level\/} of a rooted tree is the set of
vertices that are at distance $n$ from the root.  Clearly, the $n$-th level
of the rooted tree $X^*$ is formed by all words of length $n$ in the
alphabet $X$.

Now let us explain how the automaton $A$ functions.  First we choose an
{\em initial state\/} $q\in Q$ and prepare an {\em input word\/}
$w=x_1x_2\dots x_n\in X^*$.  Then we set the automaton to the state $q$ and
start inputting the word $w$ into it, letter by letter.  After reading a
letter $x'$ in a state $q'$, the automaton produces the output letter
$\psi(q',x')$ and makes transition to the state $\phi(q',x')$.  Hence the
automaton's job results in two sequences: a sequence of states
$q_0=q,q_1,\dots,q_n$, which describes the internal work of the automaton,
and the {\em output word\/} $v=y_1y_2\dots y_n\in X^*$.  Here
$q_i=\phi(q_{i-1},x_i)$ and $y_i=\psi(q_{i-1},x_i)$ for $1\le i\le n$.

For every choice of the initial state $q\in Q$ of $A$ we get a mapping
$A_q:X^*\to X^*$ that sends any input word to the corresponding output
word.  We say that $A_q$ is the transformation defined by the automaton $A$
with the initial state $q$.  Clearly, $A_q$ preserves the length of words.
Besides, $A_q$ transforms words from the left to the right, that is, the
first $n$ letters of $A_q(w)$ depend only on the first $n$ letters of $w$.
This implies that $A_q$ is an endomorphism of $X^*$ as a rooted tree.  If
$A_q$ is invertible then it belongs to the group $\Aut(X^*)$ of
automorphisms of the rooted tree $X^*$.  The set of transformations $A_q$,
$q\in Q$ is self-similar in the following sense.  For any $q\in Q$, $x\in
X$, and $w\in X^*$ we have that $A_q(xw)=yA_p(w)$, where $p=\phi(q,x)$,
$y=\psi(q,x)$.

The semigroup of transformations of $X^*$ generated by $A_q$, $q\in Q$ is
denoted by $S(A)$.  The automaton $A$ is called {\em invertible\/} if $A_q$
is invertible for all $q\in Q$.  If $A$ is invertible then $A_q$, $q\in Q$
generate a transformation group $G(A)$, which is a subgroup of $\Aut(X^*)$.
We say that $S(A)$ (resp. $G(A)$) is the semigroup (resp. group) defined by
the automaton $A$.

\begin{lemma}[\cite{VV}]\label{auto1}
Suppose the automaton $A$ is invertible.  Then the actions of the semigroup
$S(A)$ and the group $G(A)$ on $X^*$ have the same orbits.
\end{lemma}

One way to picture an automaton, which we use in this paper, is the {\em
Moore diagram}.  The Moore diagram of an automaton $A=(Q,X,\phi,\psi)$ is
a directed graph with labeled edges defined as follows.  The vertices of
the graph are states of the automaton $A$.  Every edge carries a
label of the form $x|y$, where $x,y\in X$.  The left field $x$ of the label
is referred to as the {\em input field\/} while the right field $y$ is
referred to as the {\em output field}.  The set of edges of the graph is in
a one-to-one correspondence with the set $Q\times X$.  Namely, for any
$q\in Q$ and $x\in X$ there is an edge that goes from the vertex $q$ to
$\phi(q,x)$ and carries the label $x|\psi(q,x)$.  The Moore diagram of an
automaton can have loops (edges joining a vertex to itself) and multiple
edges.  To simplify pictures, we do not draw multiple edges in this paper.
Instead, we use multiple labels.

The transformations $A_q$, $q\in Q$ can be defined in terms of the Moore
diagram of the automaton $A$.  For any $q\in Q$ and $w\in X^*$ we find a
path $\gamma$ in the Moore diagram such that $\gamma$ starts at the vertex
$q$ and the word $w$ can be obtained by reading the input fields of labels
along $\gamma$.  Such a path exists and is unique.  Then the word $A_q(w)$
is obtained by reading the output fields of labels along the path $\gamma$.

Let $\Gamma$ denote the Moore diagram of the automaton $A$.  We associate
to $\Gamma$ two directed graphs $\Gamma_1$ and $\Gamma_2$ with labeled
edges.  $\Gamma_1$ is obtained from $\Gamma$ by interchanging the input and
output fields of all labels.  That is, a label $x|y$ is replaced by $y|x$.
$\Gamma_2$ is obtained from $\Gamma$ by reversing all edges.  The {\em
inverse automaton\/} of $A$ is the automaton whose Moore diagram is
$\Gamma_1$.  The {\em reverse automaton\/} of $A$ is the automaton whose
Moore diagram is $\Gamma_2$.  The inverse and reverse automata of $A$ share
the alphabet and internal states with $A$.  Notice that any automaton is
completely determined by its Moore diagram.  However neither $\Gamma_1$ nor
$\Gamma_2$ must be the Moore diagram of an automaton.  So it is possible
that the inverse automaton or the reverse automaton (or both) of $A$ is not
well defined.

\begin{lemma}[\cite{GNS}]\label{auto2}
An automaton $A=(Q,X,\phi,\psi)$ is invertible if and only if for any $q\in
Q$ the map $\psi(q,\cdot):X\to X$ is bijective.  The inverse automaton $I$
of $A$ is well defined if and only if $A$ is invertible.  If this is the
case, then $I_q=A_q^{-1}$ for all $q\in Q$.
\end{lemma}

An automaton $A$ is called {\em reversible\/} if the reverse automaton of
$A$ is well defined.

\begin{lemma}[\cite{VV}]\label{auto3}
An automaton $A=(Q,X,\phi,\psi)$ is reversible if and only if for any $x\in
X$ the map $\phi(\cdot,x):Q\to Q$ is bijective.
\end{lemma}

Let $A=(Q,X,\phi,\psi)$ be an automaton.  For any nonempty word
$\xi=q_1q_2\dots q_n\in Q^*$ we let $A_\xi=A_{q_n}\dots A_{q_2}A_{q_1}$.
Also, we let $A_\varnothing=1$ (here $1$ stands for the unit element of the
group $\Aut(X^*)$, i.e., the identity mapping on $X^*$).  Clearly, any
element of the semigroup $S(A)$ is represented as $A_\xi$ for a nonempty
word $\xi\in Q^*$.  The map $X^*\times Q^*\to X^*$ given by $(w,\xi)\mapsto
A_\xi(w)$ defines a right action of the monoid $Q^*$ on the rooted regular
tree $X^*$.  That is, $A_{\xi_1\xi_2}(w)=A_{\xi_2}(A_{\xi_1}(w))$ for all
$\xi_1,\xi_2\in Q^*$ and $w\in X^*$.

To each finite automaton $A=(Q,X,\phi,\psi)$ we associate a {\em dual
automaton\/} $D$, which is obtained from $A$ by interchanging the alphabet
with the set of internal states and the state transition function with the
output function.  To be precise, $D=(X,Q,\tilde\phi,\tilde\psi)$, where
$\tilde\phi(x,q)=\psi(q,x)$ and $\tilde\psi(x,q)=\phi(q,x)$ for all $x\in
X$ and $q\in Q$.  Unlike the inverse and reverse automata, the dual
automaton is always well defined.  It is easy to see that $A$ is the dual
automaton of $D$.

The dual automaton $D$ defines a right action of the monoid $X^*$ on $Q^*$
given by $(\xi,w)\mapsto D_w(\xi)$.  This action and the action of $Q^*$ on
$X^*$ defined by the automaton $A$ are related in the following way.

\begin{proposition}[\cite{VV}]\label{auto4}
For any $w,u\in X^*$ and $\xi\in Q^*$,
$$
A_\xi(wu)=A_\xi(w)A_{D_w(\xi)}(u).
$$
\end{proposition}

\begin{corollary}[\cite{VV}]\label{auto5}
Suppose $A_\xi=1$ for some $\xi\in Q^*$.  Then $A_{g(\xi)}=1$ for every
$g\in S(D)$.
\end{corollary}

A finite automaton $A=(Q,X,\phi,\psi)$ is called {\em bi-reversible\/} if
the map $\phi(\cdot,x):Q\to Q$ is bijective for any $x\in X$, the map
$\psi(q,\cdot):X\to X$ is bijective for any $q\in Q$, and the map
$(\phi,\psi):Q\times X\to Q\times X$ is bijective as well.  All automata
that we consider in this paper are bi-reversible.  Below we formulate some
basic properties of bi-reversible automata (see also \cite{N}).

\begin{lemma}\label{auto6}
Given a finite automaton $A$, the following are equivalent:

(i) $A$ is bi-reversible;

(ii) $A$ is invertible, reversible, and its reverse automaton is
invertible;

(iii) $A$ is invertible, reversible, and its inverse automaton is
reversible;

(iv) $A$ is invertible, its dual automaton is invertible, and the dual
automaton of its inverse is invertible. 
\end{lemma}

\begin{proof}
Suppose $A=(Q,X,\phi,\psi)$ is a finite automaton.  By Lemma \ref{auto2},
$A$ is invertible if and only if maps $\psi(q,\cdot):X\to X$ are bijective
for all $q\in Q$.  By Lemma \ref{auto3}, $A$ is reversible if and only if
maps $\phi(\cdot,x):Q\to Q$ are bijective for all $x\in X$.  Let $\Gamma$
be the Moore diagram of $A$ and $\Gamma'$ be the graph obtained from
$\Gamma$ by reversing all edges and interchanging fields of all labels.
The graph $\Gamma'$ is the Moore diagram of an automaton if for any $q\in
Q$ and $x\in X$ there is exactly one edge of $\Gamma'$ that starts at the
vertex $q$ and has $x$ as the input field of its label.  By definition of
$\Gamma'$ the number of edges with the latter property is equal to the
number of pairs $(p,y)\in Q\times X$ such that $q=\phi(p,y)$ and
$x=\psi(p,y)$.  Therefore $\Gamma'$ is the Moore diagram of an automaton if
and only if the map $(\phi,\psi):Q\times X\to Q\times X$ is bijective.
Thus $A$ is bi-reversible if and only if it is invertible, reversible, and
$\Gamma'$ is the Moore diagram of an automaton.

Assume that the automaton $A$ is invertible and reversible.  Let $I$ and
$R$ be the inverse and reverse automata of $A$, respectively.  If the graph
$\Gamma'$ is the Moore diagram of an automaton then the automaton is both
the inverse automaton of $R$ and the reverse automaton of $I$.  On the
other hand, if $\Gamma'$ is not the Moore diagram of an automaton then $R$
is not invertible and $I$ is not reversible.  It follows that conditions
(i), (ii), and (iii) are equivalent.

It follows from Lemmas \ref{auto2} and \ref{auto3} that a finite automaton
is reversible if and only if its dual automaton is invertible.  This
implies that conditions (iii) and (iv) are equivalent.
\end{proof}

\begin{lemma}\label{auto7}
If an automaton is bi-reversible then its inverse, reverse, and dual
automata are also bi-reversible.
\end{lemma}

\begin{proof}
It follows directly from definitions that an automaton is bi-reversible if
and only if its dual automaton is bi-reversible.

Suppose $A$ is a bi-reversible automaton.  By Lemma \ref{auto6}, $A$ is
invertible and reversible.  Let $I$ and $R$ denote the inverse and reverse
automata of $A$, respectively.  By Lemma \ref{auto6}, $I$ is reversible and
$R$ is invertible.  It is easy to see that $A$ is both the inverse
automaton of $I$ and the reverse automaton of $R$.  Therefore the automata
$I$ and $R$ are invertible and reversible.  Moreover, the inverse automaton
of $I$ is reversible and the reverse automaton of $R$ is invertible.  By
Lemma \ref{auto6}, the automata $I$ and $R$ are bi-reversible.
\end{proof}

Suppose $A^{(1)}=(Q_1,X,\phi_1,\psi_1),\dots,A^{(k)}=(Q_k,X,\phi_k,\psi_k)$
are automata over the same alphabet $X$ such that their sets of internal
states $Q_1,Q_2,\dots,Q_k$ are disjoint.  The {\em disjoint union\/} of
automata $A^{(1)},A^{(2)},\dots,A^{(k)}$ is an automaton
$U=(Q_1\cup\dots\cup Q_k,X,\phi,\psi)$, where the functions $\phi$, $\psi$
are defined so that $\phi=\phi_i$ and $\psi=\psi_i$ on $Q_i\times X$ for
$1\le i\le k$.  Obviously, $U_q=A^{(i)}_q$ for all $q\in Q_i$, $1\le i\le
k$.  The Moore diagram of the automaton $U$ is the disjoint union of the
Moore diagrams of $A^{(1)},A^{(2)},\dots,A^{(k)}$.

\begin{lemma}\label{auto8}
The disjoint union of automata $A^{(1)},A^{(2)},\dots,A^{(k)}$ is
invertible (resp. reversible, bi-reversible) if and only if each $A^{(i)}$
is invertible (resp. reversible, bi-reversible).
\end{lemma}

\begin{proof}
Suppose that an automaton $U$ is the disjoint union of automata $A^{(1)},
\dots,A^{(k)}$.  Note that the disjoint union of graphs $\Gamma_1,\dots,
\Gamma_k$ is the Moore diagram of an automaton over an alphabet $X$ if and
only if each $\Gamma_i$ is the Moore diagram of an automaton defined over
$X$.  Since the Moore diagram of $U$ is the disjoint union of the Moore
diagrams of $A^{(1)},\dots,A^{(k)}$, it follows that $U$ is invertible
(resp. reversible) if and only if each $A^{(i)}$ is invertible (resp.
reversible).  Moreover, if $U$ is invertible then its inverse automaton is
the disjoint union of the inverse automata of $A^{(1)},\dots,A^{(k)}$.
Hence the inverse automaton of $U$ is reversible if and only if the inverse
automaton of each $A^{(i)}$ is reversible.  Now Lemma \ref{auto6} implies
that $U$ is bi-reversible if and only if each $A^{(i)}$ is bi-reversible.
\end{proof}

\section{The Aleshin automaton}\label{a}

In this section we recall constructions and results of the paper \cite{VV}
where the Aleshin automaton was studied.  Some constructions are slightly
modified.

The Aleshin automaton is an automaton $A$ over the alphabet $X=\{0,1\}$
with the set of internal states $Q=\{a,b,c\}$.  The state transition
function $\phi$ and the output function $\psi$ of $A$ are defined as
follows: $\phi(a,0)=\phi(b,1)=c$, $\phi(a,1)=\phi(b,0)=b$, $\phi(c,0)=
\phi(c,1)=a$; $\psi(a,0)=\psi(b,0)=\psi(c,1)=1$, $\psi(a,1)=\psi(b,1)=
\psi(c,0)=0$.  The Moore diagram of $A$ is depicted in Figure \ref{fig1}.
It is easy to verify that the automaton $A$ is invertible and reversible.
Moreover, the inverse automaton of $A$ can be obtained from $A$ by
renaming letters $0$ and $1$ of the alphabet to $1$ and $0$, respectively.
The reverse automaton of $A$ can be obtained from $A$ by renaming its
states $a$ and $c$ to $c$ and $a$, respectively.  Lemma \ref{auto6}
implies that $A$ is bi-reversible.

\begin{figure}[t]
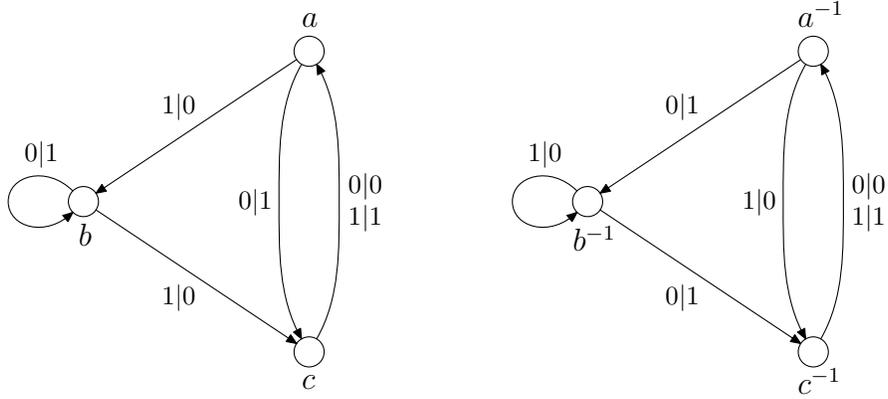

%
\includegraphics{first.eps}
\includegraphics{first.eps}
\begin{picture}(0,0)(0,80)
\put(20,0)
{
 \begin{picture}(0,0)(0,0)
 \put(119,66){$a$}
 \put(34.5,-17){$b$}
 \put(119,-71){$c$}
 \put(66,33){\footnotesize $1|0$}
 \put(66,-39){\footnotesize $1|0$}
 \put(14,15){\footnotesize $0|1$}
 \put(95,-3){\footnotesize $0|1$}
 \put(131.5,-3){\footnotesize $\begin{array}{c} 0|0 \\ 1|1 \end{array}$}
 \end{picture}
}
\put(210.5,0)
{
 \begin{picture}(0,0)(0,0)
 \put(116,66){$a^{-1}$}
 \put(31,-19){$b^{-1}$}
 \put(116,-73){$c^{-1}$}
 \put(66,33){\footnotesize $0|1$}
 \put(66,-39){\footnotesize $0|1$}
 \put(14,15){\footnotesize $1|0$}
 \put(95,-3){\footnotesize $1|0$}
 \put(131.5,-3){\footnotesize $\begin{array}{c} 0|0 \\ 1|1 \end{array}$}
 \end{picture}
}
\end{picture}
\vspace{155bp}
\caption{
\label{fig3}
Automaton $U$.
}
\end{figure}

Let $I$ denote the automaton obtained from the inverse of $A$ by renaming
its states $a$, $b$, $c$ to $a^{-1}$, $b^{-1}$, $c^{-1}$, respectively.
Here, $a^{-1}$, $b^{-1}$, and $c^{-1}$ are assumed to be elements of the
free group on generators $a$, $b$, $c$.  Further, let $U$ denote the
disjoint union of automata $A$ and $I$.  The automaton $U$ is defined over
the alphabet $X=\{0,1\}$, with the set of internal states $Q^\pm=
\{a,b,c,a^{-1},b^{-1},c^{-1}\}$.  By definition, $U_a=A_a$, $U_b=A_b$,
$U_c=A_c$, $U_{a^{-1}}=A_a^{-1}$, $U_{b^{-1}}=A_b^{-1}$, $U_{c^{-1}}=
A_c^{-1}$.

\begin{figure}[b]
%
\includegraphics{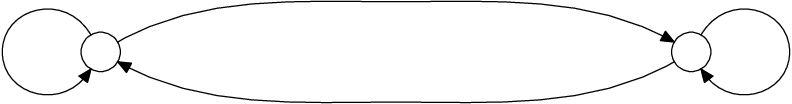}
\begin{picture}(0,0)(-103,55)
\put(-1,-16){$0$}
\put(168,-16){$1$}
\put(45,27){\small $\begin{array}{r} a|c \\ b|b \end{array}$}
\put(90,27){\small $\begin{array}{r} a^{-1}|b^{-1} \\ b^{-1}|c^{-1} \end{array}$}
\put(45,-34){\small $\begin{array}{r} a|b \\ b|c \end{array}$}
\put(90,-34){\small $\begin{array}{r} a^{-1}|c^{-1} \\ b^{-1}|b^{-1} \end{array}$}
\put(-70,-3){\small $\begin{array}{c} c|a \\ c^{-1}|a^{-1} \end{array}$}
\put(197,-3){\small $\begin{array}{c} c|a \\ c^{-1}|a^{-1} \end{array}$}
\end{picture}
\vspace{105bp}
\caption{
\label{fig4}
Automaton $D$.
}
\end{figure}

Let $D$ denote the dual automaton of the automaton $U$.  The automaton $D$
is defined over the alphabet $Q^\pm$, with two internal states $0$ and $1$.
By $\phi_D$ denote its transition function.  Then $\phi_D(0,q)=1$ and
$\phi_D(1,q)=0$ for $q\in\{a,b,a^{-1},b^{-1}\}$, while $\phi_D(0,q)=0$ and
$\phi_D(1,q)=1$ for $q\in\{c,c^{-1}\}$.  Also, we consider an auxiliary
automaton $E$ that is closely related to $D$.  By definition, the automaton
$E$ shares with $D$ the alphabet, the set of internal states, and the state
transition function.  The output function $\psi_E$ of $E$ is defined so
that $\psi_E(0,q)=\si_0(q)$ and $\psi_E(1,q)=\si_1(q)$ for all $q\in
Q^\pm$, where $\si_0=(a^{-1}b^{-1})$ and $\si_1=(ab)$ are permutations on
the set $Q^\pm$.

\begin{figure}[tb]
%
\includegraphics{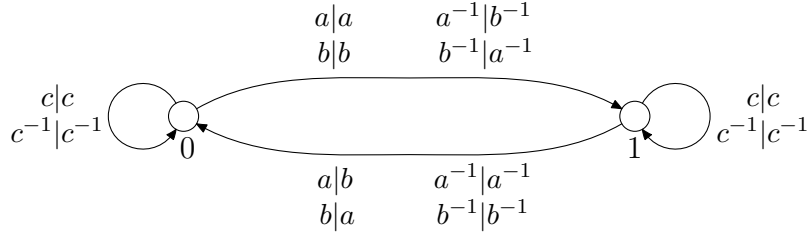}
\begin{picture}(0,0)(-103,55)
\put(-1,-16){$0$}
\put(168,-16){$1$}
\put(45,27){\small $\begin{array}{c} a|a \\ b|b \end{array}$}
\put(90,27){\small $\begin{array}{c} a^{-1}|b^{-1} \\ \,b^{-1}|a^{-1} \end{array}$}
\put(45,-34){\small $\begin{array}{c} a|b \\ \,b|a \end{array}$}
\put(90,-34){\small $\begin{array}{c} a^{-1}|a^{-1} \\ b^{-1}|b^{-1} \end{array}$}
\put(-70,-3){\small $\begin{array}{c} c|c \\ c^{-1}|c^{-1} \end{array}$}
\put(197,-3){\small $\begin{array}{c} c|c \\ c^{-1}|c^{-1} \end{array}$}
\end{picture}
\vspace{105bp}
\caption{
\label{fig5}
Automaton $E$.
}
\end{figure}

Lemmas \ref{auto7} and \ref{auto8} imply that $I$, $U$, and $D$ are
bi-reversible automata.  As for the automaton $E$, it is easy to verify
that $E$ coincides with its inverse automaton while the reverse automaton
of $E$ can be obtained from $E$ by renaming its states $0$ and $1$ to $1$
and $0$, respectively.  Hence $E$ is bi-reversible due to Lemma
\ref{auto6}.

To each permutation $\tau$ on the set $Q=\{a,b,c\}$ we assign an
automorphism $\pi_\tau$ of the free monoid $(Q^\pm)^*$.  The automorphism
$\pi_\tau$ is uniquely defined by $\pi_\tau(q)=\tau(q)$, $\pi_\tau(q^{-1})=
(\tau(q))^{-1}$ for all $q\in Q$.  Let $\langle a,b,c\rangle$ denote the
free group on generators $a$, $b$, and $c$, let $\delta:(Q^\pm)^*\to\langle
a,b,c\rangle$ be the homomorphism that sends each element of $Q^\pm\subset
(Q^\pm)^*$ to itself, and let $p_\tau$ be the automorphism of $\langle
a,b,c\rangle$ defined by $p_\tau(q)=\tau(q)$, $q\in Q$.  Then
$\delta(\pi_\tau(\xi))=p_\tau(\delta(\xi))$ for all $\xi\in(Q^\pm)^*$.

\begin{lemma}[\cite{VV}]\label{a1}
(i) $E_0^2=E_1^2=1$, $E_0E_1=E_1E_0=\pi_{(ab)}$;

(ii) $D_0=\pi_{(ac)}E_0=\pi_{(abc)}E_1$, $D_1=\pi_{(abc)}E_0=
\pi_{(ac)}E_1$.
\end{lemma}

\begin{proposition}[\cite{VV}]\label{a2}
The group $G(D)$ contains $E_0$, $E_1$, and all transformations of the form
$\pi_\tau$.  Moreover, $G(D)$ is generated by $E_0$, $\pi_{(ab)}$, and
$\pi_{(bc)}$.
\end{proposition}

As shown in Section \ref{auto}, the automaton $U$ defines a right action
$X^*\times(Q^\pm)^*\to X^*$ of the monoid $(Q^\pm)^*$ on the rooted binary
tree $X^*$ given by $(w,\xi)\mapsto U_\xi(w)$.  Let $\chi:(Q^\pm)^*\to
\{-1,1\}$ be the unique homomorphism such that $\chi(a)=\chi(b)=
\chi(a^{-1})=\chi(b^{-1})=-1$, $\chi(c)=\chi(c^{-1})=1$.

\begin{lemma}[\cite{VV}]\label{a3}
Given $\xi\in(Q^\pm)^*$, the automorphism $U_\xi$ of the rooted binary tree
$\{0,1\}^*$ acts trivially on the first level of the tree (i.e., on
one-letter words) if and only if $\chi(\xi)=1$.
\end{lemma}

Now we introduce an alphabet consisting of two symbols $*$ and $*^{-1}$.  A
word over the alphabet $\{*,*^{-1}\}$ is called a {\em pattern}.  Every
word $\xi$ over the alphabet $Q^\pm$ is assigned a pattern $v$ that is
obtained from $\xi$ by substituting $*$ for each occurrence of letters
$a,b,c$ and substituting $*^{-1}$ for each occurrence of letters
$a^{-1},b^{-1},c^{-1}$.  We say that $v$ is the pattern of $\xi$ or that
$\xi$ follows the pattern $v$.

A word $\xi=q_1q_2\dots q_n\in(Q^\pm)^*$ is called {\em freely
irreducible\/} if none of its two-letter subwords $q_1q_2,q_2q_3,\dots,
q_{n-1}q_n$ coincides with one of the following words: $aa^{-1},bb^{-1},
cc^{-1},a^{-1}a,b^{-1}b,c^{-1}c$.  Otherwise $\xi$ is called {\em freely
reducible}.

\begin{lemma}[\cite{VV}]\label{a4}
For any nonempty pattern $v$ there exist words $\xi_1,\xi_2\in(Q^\pm)^*$
such that $\xi_1$ and $\xi_2$ are freely irreducible, follow the pattern
$v$, and $\chi(\xi_2)=-\chi(\xi_1)$.
\end{lemma}

\begin{proposition}[\cite{VV}]\label{a5}
Suppose $\xi\in(Q^\pm)^*$ is a freely irreducible word.  Then the orbit of
$\xi$ under the action of the group $G(D)$ on $(Q^\pm)^*$ consists of all
freely irreducible words following the same pattern as $\xi$.
\end{proposition}

\begin{corollary}\label{a5plus}
The group defined by the dual automaton of $A$ acts transitively on each
level of the rooted ternary tree $Q^*$.
\end{corollary}

\begin{proof}
Let $D^+$ denote the dual automaton of $A$.  The rooted tree $Q^*$ is a
subtree of $(Q^\pm)^*$.  It is easy to see that $Q^*$ is invariant under
transformations $D_0$, $D_1$ and the restrictions of these transformations
to $Q^*$ are $D^+_0$, $D^+_1$.  In particular, the orbits of the $G(D^+)$
action on $Q^*$ are those orbits of the $G(D)$ action on $(Q^\pm)^*$ that
are contained in $Q^*$.  Any level of the tree $Q^*$ consists of words of a
fixed length over the alphabet $Q$.  As elements of $(Q^\pm)^*$, all these
words are freely irreducible and follow the same pattern.  Proposition
\ref{a5} implies that they are in the same orbit of the $G(D^+)$ action.
\end{proof}

Lemmas \ref{a3}, \ref{a4} and Proposition \ref{a5} lead to the following
statement.

\begin{theorem}[\cite{VV}]\label{a6}
The group $G(A)$ is the free non-Abelian group on generators $A_a$, $A_b$,
$A_c$.
\end{theorem}

\section{Series of finite automata of Aleshin type}\label{series}

In this section we consider a series of finite automata starting from the
Aleshin automaton.  We use the notation of the previous section.

For any integer $n\ge1$ we define an Aleshin type automaton $A^{(n)}$.
This is an automaton over the alphabet $X=\{0,1\}$ with a set of states
$Q_n$ of cardinality $2n+1$.  The states of $A^{(n)}$ are denoted so that
$Q_1=\{a_1,b_1,c_1\}$ and $Q_n=\{a_n,b_n,c_n,q_{n1},\dots,q_{n,2n-2}\}$ for
$n\ge2$.  The state transition function $\phi_n$ of $A^{(n)}$ is defined as
follows: $\phi_n(a_n,0)=\phi_n(b_n,1)=c_n$, $\phi_n(a_n,1)=\phi_n(b_n,0)=
b_n$, and $\phi_n(q_{ni},0)=\phi_n(q_{ni},1)=q_{n,i+1}$ for $0\le i\le
2n-2$, where by definition $q_{n0}=c_n$ and $q_{n,2n-1}=a_n$.  The output
function $\psi_n$ of $A^{(n)}$ is defined so that for any $x\in X$ we have
$\psi_n(q,x)=1-x$ if $q\in\{a_n,b_n\}$ and $\psi_n(q,x)=x$ if $q\in
Q_n\setminus\{a_n,b_n\}$.

\begin{figure}[t]
%
\includegraphics{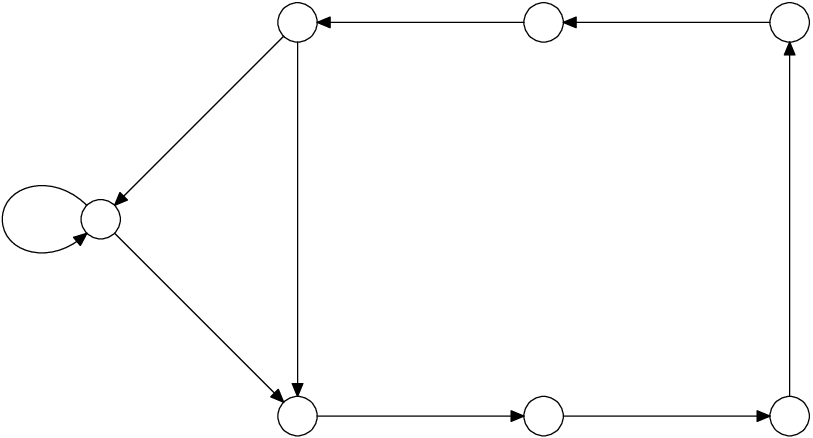}
\begin{picture}(0,0)(-54,86)
 \put(90,67){$a_3$}
 \put(34,-17){$b_3$}
 \put(90,-71){$c_3$}
 \put(160,-71){$q_{31}$}
 \put(230,-71){$q_{32}$}
 \put(230,67.5){$q_{33}$}
 \put(160,67.5){$q_{34}$}
 \put(53,33){\footnotesize $1|0$}
 \put(53,-39){\footnotesize $1|0$}
 \put(14,15){\footnotesize $0|1$}
 \put(77.5,-3){\footnotesize $0|1$}
 \put(214,-3){\footnotesize $\begin{array}{c} 0|0 \\ 1|1 \end{array}$}
 \put(118,-45.5){\footnotesize $\begin{array}{c} 0|0 \\ 1|1 \end{array}$}
 \put(118,39.5){\footnotesize $\begin{array}{c} 0|0 \\ 1|1 \end{array}$}
 \put(189,-45.5){\footnotesize $\begin{array}{c} 0|0 \\ 1|1 \end{array}$}
 \put(189,39.5){\footnotesize $\begin{array}{c} 0|0 \\ 1|1 \end{array}$}
\end{picture}
\vspace{162bp}
\caption{
\label{fig6}
Automaton $A^{(3)}$.
}
\end{figure}

Up to renaming of the internal states, $A^{(1)}$ and $A^{(2)}$ are the two
automata introduced by Aleshin \cite{A} (see Figure \ref{fig1}).

We shall deal with automata $A^{(n)}$ by following the framework developed
in the paper \cite{VV} and described in Section \ref{a}.

Let us fix a positive integer $n$.  It is easy to see that the inverse
automaton of the automaton $A^{(n)}$ can be obtained from $A^{(n)}$ by
renaming letters $0$ and $1$ of the alphabet to $1$ and $0$, respectively.
Besides, the reverse automaton of $A^{(n)}$ can be obtained from $A^{(n)}$
by renaming its states $c_n,q_{n1},\dots,q_{n,2n-2},a_n$ to
$a_n,q_{n,2n-2},\dots,q_{n1},c_n$, respectively.  Lemma \ref{auto6} implies
that $A^{(n)}$ is bi-reversible.

Let $I^{(n)}$ denote the automaton obtained from the inverse of $A^{(n)}$
by renaming each state $q\in Q_n$ to $q^{-1}$, where $q^{-1}$ is regarded
as an element of the free group on generators $a_n,b_n,c_n,q_{n1},\dots,
q_{n,2n-2}$.  Further, let $U^{(n)}$ denote the disjoint union of automata
$A^{(n)}$ and $I^{(n)}$.  The automaton $U^{(n)}$ is defined over the
alphabet $X=\{0,1\}$, with the set of internal states $Q_n^\pm=
\bigcup_{q\in Q_n}\{q,q^{-1}\}$.  By definition, $U^{(n)}_q=A^{(n)}_q$ and
$U^{(n)}_{q^{-1}}=(A^{(n)}_q)^{-1}$ for all $q\in Q_n$.

Let $D^{(n)}$ denote the dual automaton of the automaton $U^{(n)}$.  The
automaton $D^{(n)}$ is defined over the alphabet $Q_n^\pm$, with two
internal states $0$ and $1$.  By $\la_n$ denote its transition function.
Then $\la_n(0,q)=1$ and $\la_n(1,q)=0$ if $q\in\{a_n,b_n,a_n^{-1},
b_n^{-1}\}$ while $\la_n(0,q)=0$ and $\la_n(1,q)=1$ otherwise.  Also, we
consider an auxiliary automaton $E^{(n)}$.  By definition, the automaton
$E^{(n)}$ shares with $D^{(n)}$ the alphabet, the set of internal states,
and the state transition function.  The output function $\mu_n$ of
$E^{(n)}$ is defined so that $\mu_n(0,q)=\si_0(q)$ and $\mu_n(1,q)=
\si_1(q)$ for all $q\in Q_n^\pm$, where $\si_0=(a_n^{-1}b_n^{-1})$ and
$\si_1=(a_nb_n)$ are permutations on the set $Q_n^\pm$.

Lemmas \ref{auto7} and \ref{auto8} imply that $I^{(n)}$, $U^{(n)}$, and
$D^{(n)}$ are bi-reversible automata.  Further, it is easy to see that the
automaton $E^{(n)}$ coincides with its inverse automaton while the reverse
automaton of $E^{(n)}$ can be obtained from $E^{(n)}$ by renaming its
states $0$ and $1$ to $1$ and $0$, respectively.  By Lemma \ref{auto6},
$E^{(n)}$ is bi-reversible.

To each permutation $\tau$ on the set $Q_n$ we assign an automorphism
$\pi^{(n)}_\tau$ of the free monoid $(Q_n^\pm)^*$ such that
$\pi^{(n)}_\tau(q)=\tau(q)$, $\pi^{(n)}_\tau(q^{-1})=(\tau(q))^{-1}$ for
all $q\in Q_n$.  The automorphism $\pi^{(n)}_\tau$ is uniquely determined
by $\tau$.

\begin{lemma}\label{series1}
(i) $(E^{(n)}_0)^2=(E^{(n)}_1)^2=1$, $E^{(n)}_0E^{(n)}_1=
E^{(n)}_1E^{(n)}_0=\pi^{(n)}_{(a_nb_n)}$;

(ii) $D^{(n)}_0=\pi^{(n)}_{\tau_0}E^{(n)}_0=\pi^{(n)}_{\tau_1}E^{(n)}_1$,
$D^{(n)}_1=\pi^{(n)}_{\tau_1}E^{(n)}_0=\pi^{(n)}_{\tau_0}E^{(n)}_1$, where
$\tau_0=(a_nc_nq_{n1}\dots q_{n,2n-2})$, $\tau_1=(a_nb_nc_nq_{n1}\dots
q_{n,2n-2})$.
\end{lemma}

\begin{proof}
Since the inverse automaton of $E^{(n)}$ coincides with $E^{(n)}$, Lemma
\ref{auto2} implies that $(E^{(n)}_0)^2=(E^{(n)}_1)^2=1$.

We have that $E^{(n)}=(X,Q_n^\pm,\la_n,\mu_n)$, where the functions $\la_n$
and $\mu_n$ are defined above.  Note that the function $\la_n$ does not
change when elements $0$ and $1$ of the set $X$ are renamed to $1$ and $0$,
respectively.  For any permutation $\si$ on the set $Q_n^\pm$ we define an
automaton $Y^\si=(X,Q_n^\pm,\la_n,\si\mu_n)$.  The Moore diagram of $Y^\si$
is obtained from the Moore diagram of $E^{(n)}$ by applying $\si$ to the
output fields of all labels.  It is easy to observe that $Y^\tau_0=
\alpha_\si E^{(n)}_0$ and $Y^\tau_1=\alpha_\si E^{(n)}_1$, where
$\alpha_\si$ is the unique automorphism of the monoid $(Q_n^\pm)^*$ such
that $\alpha_\si(q)=\si(q)$ for all $q\in Q_n^\pm$.

Let us consider the following permutations on $Q_n^\pm$:
\begin{eqnarray*}
& \si_0=(a_n^{-1}b_n^{-1}), \qquad \si_1=(a_nb_n), \qquad
\si_2=(a_nb_n)(a_n^{-1}b_n^{-1}),\\
& \si_3=(a_nc_nq_{n1}\dots q_{n,2n-2})(a_n^{-1}b_n^{-1}c_n^{-1}q_{n1}^{-1}
\dots q_{n,2n-2}^{-1}),\\
& \si_4=(a_nb_nc_nq_{n1}\dots q_{n,2n-2})(a_n^{-1}c_n^{-1}q_{n1}^{-1}\dots
q_{n,2n-2}^{-1}),\\
& \si_5=(a_nc_nq_{n1}\dots q_{n,2n-2})(a_n^{-1}c_n^{-1}q_{n1}^{-1}\dots
q_{n,2n-2}^{-1}),\\
& \si_6=(a_nb_nc_nq_{n1}\dots q_{n,2n-2})
(a_n^{-1}b_n^{-1}c_n^{-1}q_{n1}^{-1}\dots q_{n,2n-2}^{-1}).
\end{eqnarray*}
Since $\si_2\si_0=\si_1$ and $\si_2\si_1=\si_0$, it follows that the
automaton $Y^{\si_2}$ can be obtained from $E^{(n)}$ by renaming its states
$0$ and $1$ to $1$ and $0$, respectively.  Therefore $E^{(n)}_0=
Y^{\si_2}_1=\alpha_{\si_2}E^{(n)}_1$ and $E^{(n)}_1=Y^{\si_2}_0=
\alpha_{\si_2}E^{(n)}_0$.  Consequently, $E^{(n)}_0E^{(n)}_1=\alpha_{\si_2}
(E^{(n)}_1)^2=\alpha_{\si_2}$ and $E^{(n)}_1E^{(n)}_0=\alpha_{\si_2}
(E^{(n)}_0)^2=\alpha_{\si_2}$.  Clearly, $\alpha_{\si_2}=
\pi^{(n)}_{(a_nb_n)}$.

Since $\si_5\si_0=\si_3$ and $\si_5\si_1=\si_4$, it follows that
$Y^{\si_5}=D^{(n)}$.  Hence $D^{(n)}_0=\alpha_{\si_5}E^{(n)}_0$ and
$D^{(n)}_1=\alpha_{\si_5}E^{(n)}_1$.  Furthermore, the equalities
$\si_6\si_0=\si_4$ and $\si_6\si_1=\si_3$ imply that the automaton
$Y^{\si_6}$ can be obtained from $D^{(n)}$ by renaming its states $0$ and
$1$ to $1$ and $0$, respectively.  Therefore $D^{(n)}_0=Y^{\si_6}_1=
\alpha_{\si_6}E^{(n)}_1$ and $D^{(n)}_1=Y^{\si_6}_0=\alpha_{\si_6}
E^{(n)}_0$.  It remains to notice that $\alpha_{\si_5}=\pi^{(n)}_{\tau_0}$
and $\alpha_{\si_6}=\pi^{(n)}_{\tau_1}$.
\end{proof}

\begin{lemma}\label{series2}
For any integer $M\ge3$ the group of permutations on the set $\{1,2,\dots,
M\}$ is generated by permutations $(12)$ and $(123\dots M)$.
\end{lemma}

\begin{proof}
Let $\tau_0=(12)$, $\tau_1=(123\dots M)$, and $\tau_2=(23\dots M)$.  Then
$\tau_2=\tau_0\tau_1$.  For any $k$, $2\le k\le M$ we have $(1k)=
\tau_2^{k-2}\tau_0\tau_2^{-(k-2)}$.  Further, for any $l$ and $m$, $1\le l<
m\le M$ we have $(lm)=\tau_1^{l-1}(1k)\tau_1^{-(l-1)}$, where $k=m-l+1$.
Therefore the group generated by $\tau_0$ and $\tau_1$ contains all
transpositions $(lm)$, $1\le l<m\le M$.  It remains to notice that any
permutation on $\{1,2,\dots,M\}$ is a product of transpositions.
\end{proof}

\begin{proposition}\label{series3}
The group $G(D^{(n)})$ contains $E^{(n)}_0$, $E^{(n)}_1$, and all
transformations of the form $\pi^{(n)}_\tau$.  Moreover, $G(D^{(n)})$ is
generated by $E^{(n)}_0$, $\pi^{(n)}_{\tau_0}$, and $\pi^{(n)}_{\tau_1}$,
where $\tau_0=(a_nc_nq_{n1}\dots q_{n,2n-2})$, $\tau_1=
(a_nb_nc_nq_{n1}\dots q_{n,2n-2})$.
\end{proposition}

\begin{proof}
It is easy to see that $\pi^{(n)}_{\tau\si}=\pi^{(n)}_\tau\pi^{(n)}_\si$
for any permutations $\tau$ and $\si$ on the set $Q_n$.  It follows that
$\pi^{(n)}_{\tau^{-1}}=(\pi^{(n)}_\tau)^{-1}$ for any permutation $\tau$ on
$Q_n$.

By Lemma \ref{series1}, the group generated by $E^{(n)}_0$,
$\pi^{(n)}_{\tau_0}$, and $\pi^{(n)}_{\tau_1}$ contains $G(D^{(n)})$.
Besides, $D^{(n)}_0(D^{(n)}_1)^{-1}=\pi^{(n)}_{\tau_0}E^{(n)}_0
(\pi^{(n)}_{\tau_1}E^{(n)}_0)^{-1}=\pi^{(n)}_{\tau_0}
(\pi^{(n)}_{\tau_1})^{-1}$.  By the above $\pi^{(n)}_{\tau_0}
(\pi^{(n)}_{\tau_1})^{-1}=\pi^{(n)}_{\tau_2}$, where $\tau_2=
\tau_0\tau_1^{-1}=(b_nc_n)$.  Similarly,
$$
(D^{(n)}_0)^{-1}D^{(n)}_1=(\pi^{(n)}_{\tau_0}E^{(n)}_0)^{-1}
\pi^{(n)}_{\tau_1}E^{(n)}_0=(E^{(n)}_0)^{-1}\pi^{(n)}_{\tau_3}E^{(n)}_0,
$$
where $\tau_3=\tau_0^{-1}\tau_1=(a_nb_n)$.  Lemma \ref{series1} implies
that $E^{(n)}_0$ and $\pi^{(n)}_{\tau_3}$ commute, hence
$(D^{(n)}_0)^{-1}D^{(n)}_1=\pi^{(n)}_{\tau_3}$.  Consider two more
permutations on $Q_n$: $\tau_4=(a_nc_n)$ and $\tau_5=(c_nq_{n1}\dots
q_{n,2n-2})$.  Note that $\tau_4=\tau_2\tau_3\tau_2$ and $\tau_5=
\tau_4\tau_0$.  By the above $\pi^{(n)}_{\tau_2},\pi^{(n)}_{\tau_3}\in
G(D^{(n)})$, hence $\pi^{(n)}_{\tau_4}\in G(D^{(n)})$.  Then
$\pi^{(n)}_{\tau_5}E^{(n)}_0=\pi^{(n)}_{\tau_4}\pi^{(n)}_{\tau_0}
E^{(n)}_0=\pi^{(n)}_{\tau_4}D^{(n)}_0\in G(D^{(n)})$.  Since $\tau_5(a_n)=
a_n$ and $\tau_5(b_n)=b_n$, it easily follows that transformations
$\pi^{(n)}_{\tau_5}$ and $E^{(n)}_0$ commute.  As $\tau_5$ is a permutation
of odd order $2n-1$ while $E^{(n)}_0$ is an involution, we have that
$(\pi^{(n)}_{\tau_5}E^{(n)}_0)^{2n-1}=E^{(n)}_0$.  In particular,
$E^{(n)}_0\in G(D^{(n)})$.  Now Lemma \ref{series1} implies that
$\pi^{(n)}_{\tau_0},\pi^{(n)}_{\tau_1},E^{(n)}_1\in G(D^{(n)})$.

By Lemma \ref{series2}, the group of all permutations on the set $Q_n$ is
generated by permutations $\tau_1$ and $\tau_3$.  Since
$\pi^{(n)}_{\tau_1},\pi^{(n)}_{\tau_3}\in G(D^{(n)})$, it follows that
$G(D^{(n)})$ contains all transformations of the form $\pi^{(n)}_\tau$.
\end{proof}

Recall that words over the alphabet $\{*,*^{-1}\}$ are called patterns.
Every word $\xi\in(Q_n^\pm)^*$ is assigned a pattern $v$ that is obtained
from $\xi$ by substituting $*$ for each occurrence of letters $a_n,b_n,c_n,
q_{n1},\dots,q_{n,2n-2}$ and substituting $*^{-1}$ for each occurrence of
letters $a_n^{-1},b_n^{-1},c_n^{-1},q_{n1}^{-1},\dots,q_{n,2n-2}^{-1}$.  We
say that $\xi$ follows the pattern $v$.

A word $\xi=q_1q_2\dots q_k\in(Q_n^\pm)^*$ is called freely irreducible if
none of its two-letter subwords $q_1q_2,q_2q_3,\dots,q_{k-1}q_k$ is of the
form $qq^{-1}$ or $q^{-1}q$, where $q\in Q_n$.  Otherwise $\xi$ is called
freely reducible.

\begin{lemma}\label{series4}
For any nonempty pattern $v$ there exists a freely irreducible word
$\xi\in(Q_n^\pm)^*$ such that $v$ is the pattern of $\xi$ and the
transformation $U^{(n)}_\xi$ acts nontrivially on the first level of the
rooted binary tree $X^*$.
\end{lemma}

\begin{proof}
Given a nonempty pattern $v$, let us substitute $a_n$ for each occurrence
of $*$ in $v$ and $b_n^{-1}$ for each occurrence of $*^{-1}$.  We get a
word $\xi\in(Q_n^\pm)^*$ that follows the pattern $v$.  Now let us modify
$\xi$ by changing its last letter.  If this letter is $a_n$, we change it
to $c_n$.  If the last letter of $\xi$ is $b_n^{-1}$, we change it to
$c_n^{-1}$.  This yields another word $\eta\in(Q_n^\pm)^*$ that follows the
pattern $v$.  By construction, $\xi$ and $\eta$ are freely irreducible.
Furthermore, $U^{(n)}_\eta=A^{(n)}_{c_n}(A^{(n)}_{a_n})^{-1}U^{(n)}_\xi$ if
the last letter of $v$ is $*$ while $U^{(n)}_\eta=(A^{(n)}_{c_n})^{-1}
A^{(n)}_{b_n}U^{(n)}_\xi$ if the last letter of $v$ is $*^{-1}$.  Both
$A^{(n)}_{c_n}(A^{(n)}_{a_n})^{-1}$ and $(A^{(n)}_{c_n})^{-1}A^{(n)}_{b_n}$
interchange one-letter words $0$ and $1$.  It follows that one of the
transformations $U^{(n)}_\xi$ and $U^{(n)}_\eta$ also acts nontrivially on
the first level of the rooted tree $\{0,1\}^*$.
\end{proof}

Given a nonempty, freely irreducible word $\xi\in(Q_n^\pm)^*$, let
$Z_n(\xi)$ denote the set of all freely irreducible words in $(Q_n^\pm)^*$
that follow the same pattern as $\xi$ and match $\xi$ completely or except
for the last letter.  Obviously, $\xi\in Z_n(\xi)$, and $\eta\in Z_n(\xi)$
if and only if $\xi\in Z_n(\eta)$.  The set $Z_n(\xi)$ consists of $2n$ or
$2n+1$ words.  Namely, there are exactly $2n+1$ words in $(Q_n^\pm)^*$ that
follow the same pattern as $\xi$ and match $\xi$ completely or except for
the last letter.  However if the last two letters in the pattern of $\xi$
are distinct then one of these $2n+1$ words is freely reducible.

\begin{lemma}\label{series5}
For any nonempty pattern $v$ there exists a freely irreducible word $\xi\in
(Q_n^\pm)^*$ such that $v$ is the pattern of $\xi$ and the set $Z_n(\xi)$
is contained in one orbit of the $G(D^{(n)})$ action on $(Q_n^\pm)^*$.
\end{lemma}

\begin{proof}
Let $h_n:(Q^\pm)^*\to(Q_n^\pm)^*$ be the homomorphism of monoids such that
$h_n(a)=a_n$, $h_n(b)=b_n$, $h_n(c)=c_n$, $h_n(a^{-1})=a_n^{-1}$,
$h_n(b^{-1})=b_n^{-1}$, $h_n(c^{-1})=c_n^{-1}$.  The range of $h_n$
consists of words over alphabet $\{a_n,b_n,c_n,a_n^{-1},b_n^{-1},
c_n^{-1}\}$.  For any $\zeta\in(Q^\pm)^*$ the word $h_n(\zeta)$ follows the
same pattern as $\zeta$.  Besides, $h_n(\zeta)$ is freely irreducible if
and only if $\zeta$ is.  It is easy to see that $h_n(\pi_{(ab)}(\zeta))=
\pi^{(n)}_{(a_nb_n)}(h_n(\zeta))$, $h_n(\pi_{(bc)}(\zeta))=
\pi^{(n)}_{(b_nc_n)}(h_n(\zeta))$, and $h_n(E_0(\zeta))=
E^{(n)}_0(h_n(\zeta))$.  By Proposition \ref{a2}, the group $G(D)$ is
generated by $\pi_{(ab)}$, $\pi_{(bc)}$, and $E_0$.  On the other hand,
$\pi^{(n)}_{(a_nb_n)},\pi^{(n)}_{(b_nc_n)},E^{(n)}_0\in G(D^{(n)})$ due to
Proposition \ref{series3}.  It follows that for any $g_0\in G(D)$ there
exists $g\in G(D^{(n)})$ such that $h_n(g_0(\zeta))=g(h_n(\zeta))$ for all
$\zeta\in(Q^\pm)^*$.  Now Proposition \ref{a5} implies that two words over
alphabet $\{a_n,b_n,c_n,a_n^{-1},b_n^{-1},c_n^{-1}\}$ are in the same orbit
of the $G(D^{(n)})$ action on $(Q_n^\pm)^*$ whenever they are freely
irreducible and follow the same pattern.

Let $v_0$ be the pattern obtained by deleting the last letter of $v$.  We
substitute $a_n$ for each occurrence of $*$ in $v_0$ and $b_n^{-1}$ for
each occurrence of $*^{-1}$.  This yields a word $\eta\in(Q_n^\pm)^*$ that
follows the pattern $v_0$.  Now let $\xi=\eta c_n$ if the last letter of
$v$ is $*$ and let $\xi=\eta c_n^{-1}$ otherwise.  Clearly, $\xi$ is a
freely irreducible word following the pattern $v$.  Take any $\zeta\in
Z_n(\xi)$.  If both $\zeta$ and $\xi$ are words over alphabet $\{a_n,b_n,
c_n,a_n^{-1},b_n^{-1},c_n^{-1}\}$, then it follows from the above that
$\zeta=g(\xi)$ for some $g\in G(D^{(n)})$.  Otherwise the last letter of
$\zeta$ is $q_{ni}$ or $q_{ni}^{-1}$, where $1\le i\le 2n-2$.  In this case
we have $\zeta=(\pi^{(n)}_\tau)^i(\xi)$, where $\tau=(c_nq_{n1}\dots
q_{n,2n-2})$.  By Proposition \ref{series3}, $\pi^{(n)}_\tau\in
G(D^{(n)})$.
\end{proof}

\begin{proposition}\label{series6}
Suppose $\xi\in(Q_n^\pm)^*$ is a freely irreducible word.  Then the orbit
of $\xi$ under the action of the group $G(D^{(n)})$ on $(Q_n^\pm)^*$
consists of all freely irreducible words following the same pattern as
$\xi$.
\end{proposition}

\begin{proof}
First we shall show that the $G(D^{(n)})$ action on $(Q_n^\pm)^*$
preserves patterns and free irreducibility of words.  Let $\phi_n^\pm$ and
$\psi_n^\pm$ denote the state transition and output functions of the
automaton $U^{(n)}$.  By $\tilde\phi_n$ and $\tilde\psi_n$ denote the state
transition and output functions of its dual $D^{(n)}$.  Take any $q\in
Q_n^\pm$ and $x\in X$.  By definition of $U^{(n)}$ we have that
$\phi_n^\pm(q,x)\in Q_n$ if and only if $q\in Q_n$.  Since
$\phi_n^\pm(q,x)=\tilde\psi_n(x,q)$, it follows that transformations
$D^{(n)}_0$ and $D^{(n)}_1$ preserve patterns of words.  So does any $g\in
G(D^{(n)})$.  Further, let $p=\phi_n^\pm(q,x)$ and $y=\psi_n^\pm(q,x)$.
Then $\phi_n^\pm(q^{-1},y)=p^{-1}$ and $\psi_n^\pm(q^{-1},y)=x$.
Consequently, $D^{(n)}_x(qq^{-1})=\tilde\psi_n(x,q)
\tilde\psi_n(\tilde\phi_n(x,q),q^{-1})=\phi_n^\pm(q,x)
\phi_n^\pm(q^{-1},\psi_n^\pm(q,x))=pp^{-1}$.  It follows that the set
$P=\{qq^{-1}\mid q\in Q_n^\pm\}\subset(Q_n^\pm)^*$ is invariant under
$D^{(n)}_0$ and $D^{(n)}_1$.  Any freely reducible word $\xi\in(Q_n^\pm)^*$
is represented as $\xi_1\xi_0\xi_2$, where $\xi_0\in P$ and $\xi_1,\xi_2\in
(Q_n^\pm)^*$.  For any $x\in X$ we have $D^{(n)}_x(\xi)=D^{(n)}_x(\xi_1)
D^{(n)}_{x_0}(\xi_0)D^{(n)}_{x_1}(\xi_2)$, where $x_0,x_1\in X$.  By the
above $D^{(n)}_x(\xi)$ is freely reducible.  Thus $D^{(n)}_0$ and
$D^{(n)}_1$ preserve free reducibility of words.  Since these
transformations are invertible, they also preserve free irreducibility, and
so does any $g\in G(D^{(n)})$.

Now we are going to prove that for any freely irreducible words
$\xi_1,\xi_2\in(Q_n^\pm)^*$ following the same pattern $v$ there exists
$g\in G(D^{(n)})$ such that $\xi_2=g(\xi_1)$.  The claim is proved by
induction on the length of the pattern $v$.  The empty pattern is followed
only by the empty word.  Now let $k\ge1$ and assume that the claim holds
for all patterns of length less than $k$.  Take any pattern $v$ of length
$k$.  By Lemma \ref{series5}, the pattern $v$ is followed by a freely
irreducible word $\xi\in(Q_n^\pm)^*$ such that the set $Z_n(\xi)$ is
contained in an orbit of the $G(D^{(n)})$ action.  Suppose $\xi_1,\xi_2\in
(Q_n^\pm)^*$ are freely irreducible words following the pattern $v$.  Let
$\eta,\eta_1,\eta_2$ be the words obtained by deleting the last letter of
$\xi,\xi_1,\xi_2$, respectively.  Then $\eta,\eta_1,\eta_2$ are freely
irreducible and follow the same pattern of length $k-1$.  By the inductive
assumption there are $g_1,g_2\in G(D^{(n)})$ such that $\eta=g_1(\eta_1)=
g_2(\eta_2)$.  Since the $G(D^{(n)})$ action preserves patterns and free
irreducibility, it follows that $g_1(\xi_1),g_2(\xi_2)\in Z_n(\xi)$.  As
$Z_n(\xi)$ is contained in an orbit, there exists $g_0\in G(D^{(n)})$ such
that $g_0(g_1(\xi_1))=g_2(\xi_2)$.  Then $\xi_2=g(\xi_1)$, where
$g=g_2^{-1}g_0g_1\in G(D^{(n)})$.
\end{proof}

\begin{corollary}\label{series6plus}
The group defined by the dual automaton of $A^{(n)}$ acts transitively on
each level of the rooted tree $Q_n^*$.
\end{corollary}

Corollary \ref{series6plus} follows from Proposition \ref{series6} in the
same way as Corollary \ref{a5plus} follows from Proposition \ref{a5}.  We
omit the proof.

\begin{theorem}\label{series7}
The group $G(A^{(n)})$ is the free non-Abelian group on $2n+1$ generators
$A^{(n)}_q$, $q\in Q_n$.
\end{theorem}

\begin{proof}
The group $G(A^{(n)})$ is the free non-Abelian group on generators $A_q$,
$q\in Q_n$ if and only if $(A^{(n)}_{q_1})^{m_1}(A^{(n)}_{q_2})^{m_2}\dots
(A^{(n)}_{q_k})^{m_k}\ne1$ for any pair of sequences $q_1,\dots,q_k$ and
$m_1,\dots,m_k$ such that $k>0$, $q_i\in Q_n$ and $m_i\in\bZ\setminus\{0\}$
for $1\le i\le k$, and $q_i\ne q_{i+1}$ for $1\le i\le k-1$.  Since
$U^{(n)}_q=A^{(n)}_q$ and $U^{(n)}_{q^{-1}}=(A^{(n)}_q)^{-1}$ for all $q\in
Q_n$, an equivalent condition is that $U^{(n)}_\xi\ne1$ for any nonempty
freely irreducible word $\xi\in(Q_n^\pm)^*$.

Suppose $U^{(n)}_\xi=1$ for some freely irreducible word $\xi\in
(Q_n^\pm)^*$.  By Corollary \ref{auto5}, $U^{(n)}_{g(\xi)}=1$ for all $g\in
S(D^{(n)})$.  Then Proposition \ref{auto1} imply that $U^{(n)}_{g(\xi)}=1$
for all $g\in G(D^{(n)})$.  Now it follows from Proposition \ref{series6}
that $U^{(n)}_\eta=1$ for any freely irreducible word $\eta\in(Q_n^\pm)^*$
following the same pattern as $\xi$.  In particular, $U^{(n)}_\eta$ acts
trivially on the first level of the rooted binary tree $\{0,1\}^*$.
Finally, Lemma \ref{series4} implies that $\xi$ follows the empty pattern.
Then $\xi$ itself is the empty word.
\end{proof}

\section{Disjoint unions}\label{union}

In this section we consider disjoint unions of Aleshin type automata.  We
use the notation of Sections \ref{a} and \ref{series}.

Let $N$ be a nonempty set of positive integers.  We denote by $A^{(N)}$ the
disjoint union of automata $A^{(n)}$, $n\in N$.  Then $A^{(N)}$ is an
automaton over the alphabet $X=\{0,1\}$ with the set of internal states
$Q_N=\bigcup_{n\in N}Q_n$.  It is bi-reversible since each $A^{(n)}$ is
bi-reversible.

Let $I^{(N)}$ denote the disjoint union of automata $I^{(n)}$, $n\in N$.
The automaton $I^{(N)}$ can be obtained from the inverse of $A^{(N)}$ by
renaming each state $q\in Q_N$ to $q^{-1}$.  Further, let $U^{(N)}$ denote
the disjoint union of automata $A^{(N)}$ and $I^{(N)}$.  Obviously, the
automaton $U^{(N)}$ is the disjoint union of automata $U^{(n)}$, $n\in N$.
It is defined over the alphabet $X=\{0,1\}$, with the set of internal
states $Q_N^\pm=\bigcup_{n\in N}Q_n^\pm$.  Clearly, $U^{(N)}_q=A^{(N)}_q$
and $U^{(N)}_{q^{-1}}=(A^{(N)}_q)^{-1}$ for all $q\in Q_N$.

Let $D^{(N)}$ denote the dual automaton of the automaton $U^{(N)}$.  The
automaton $D^{(N)}$ is defined over the alphabet $Q_N^\pm$, with two
internal states $0$ and $1$.  Also, we consider an auxiliary automaton
$E^{(N)}$.  By definition, the automaton $E^{(N)}$ shares with $D^{(N)}$
the alphabet, the set of internal states, and the state transition
function.  The output function $\mu_N$ of $E^{(N)}$ is defined so that
$\mu_N(0,q)=\si_0(q)$ and $\mu_N(1,q)=\si_1(q)$ for all $q\in Q_N^\pm$,
where $\si_0=\prod_{n\in N}(a_n^{-1}b_n^{-1})$ and $\si_1=\prod_{n\in N}
(a_nb_n)$ are permutations on the set $Q_N^\pm$.

Lemmas \ref{auto7} and \ref{auto8} imply that $I^{(N)}$, $U^{(N)}$, and
$D^{(N)}$ are bi-reversible automata.  Further, it is easy to see that the
automaton $E^{(N)}$ coincides with its inverse automaton while the reverse
automaton of $E^{(N)}$ can be obtained from $E^{(N)}$ by renaming its
states $0$ and $1$ to $1$ and $0$, respectively.  By Lemma \ref{auto6},
$E^{(N)}$ is bi-reversible.

To each permutation $\tau$ on the set $Q_N$ we assign an automorphism
$\pi^{(N)}_\tau$ of the free monoid $(Q_N^\pm)^*$ such that
$\pi^{(N)}_\tau(q)=\tau(q)$, $\pi^{(N)}_\tau(q^{-1})=(\tau(q))^{-1}$ for
all $q\in Q_N$.  The automorphism $\pi^{(N)}_\tau$ is uniquely determined
by $\tau$.

\begin{lemma}\label{union1}
(i) $(E^{(N)}_0)^2=(E^{(N)}_1)^2=1$, $E^{(N)}_0E^{(N)}_1=E^{(N)}_1
E^{(N)}_0=\pi^{(N)}_\tau$, where $\tau=\prod_{n\in N}(a_nb_n)$;

(ii) $D^{(N)}_0=\pi^{(N)}_{\tau_0}E^{(N)}_0=\pi^{(N)}_{\tau_1}E^{(N)}_1$,
$D^{(N)}_1=\pi^{(N)}_{\tau_1}E^{(N)}_0=\pi^{(N)}_{\tau_0}E^{(N)}_1$, where
$\tau_0=\prod_{n\in N}(a_nc_nq_{n1}\dots q_{n,2n-2})$, $\tau_1=\prod_{n\in
N}(a_nb_nc_nq_{n1}\dots q_{n,2n-2})$.
\end{lemma}

The proof of Lemma \ref{union1} is completely analogous to the proof of
Lemma \ref{series1} and we omit it.

\begin{proposition}\label{union2}
The group $G(D^{(N)})$ contains transformations $E^{(N)}_0$, $E^{(N)}_1$,
$\pi^{(N)}_{\tau_1}$, $\pi^{(N)}_{\tau_2}$, $\pi^{(N)}_{\tau_3}$, and
$\pi^{(N)}_{\tau_4}$, where $\tau_1=\prod_{n\in N}(a_nb_nc_nq_{n1}\dots
q_{n,2n-2})$, $\tau_2=\prod_{n\in N}(c_nq_{n1}\dots q_{n,2n-2})$, $\tau_3=
\prod_{n\in N}(a_nb_n)$, and $\tau_4=\prod_{n\in N}(b_nc_n)$.
\end{proposition}

\begin{proof}
It is easy to see that $\pi^{(N)}_{\tau\si}=\pi^{(N)}_\tau\pi^{(N)}_\si$
for any permutations $\tau$ and $\si$ on the set $Q_N$.  It follows that
$\pi^{(N)}_{\tau^{-1}}=(\pi^{(N)}_\tau)^{-1}$ for any permutation $\tau$ on
$Q_N$.

By Lemma \ref{union1}, $D^{(N)}_0(D^{(N)}_1)^{-1}=\pi^{(N)}_{\tau_0}
E^{(N)}_0(\pi^{(N)}_{\tau_1}E^{(N)}_0)^{-1}=\pi^{(N)}_{\tau_0}
(\pi^{(N)}_{\tau_1})^{-1}$, where $\tau_0=\prod_{n\in N}(a_nc_nq_{n1}\dots
q_{n,2n-2})$.  Since $\tau_0\tau_1^{-1}=\tau_4$, it follows that $D^{(N)}_0
(D^{(N)}_1)^{-1}=\pi^{(N)}_{\tau_4}$.  Similarly,
$$
(D^{(N)}_0)^{-1}D^{(N)}_1=(\pi^{(N)}_{\tau_0}E^{(N)}_0)^{-1}
\pi^{(N)}_{\tau_1}E^{(N)}_0=(E^{(N)}_0)^{-1}\pi^{(N)}_{\tau_3}E^{(N)}_0
$$
since $\tau_3=\tau_0^{-1}\tau_1$.  Lemma \ref{union1} implies that
$E^{(N)}_0$ and $\pi^{(N)}_{\tau_3}$ commute, hence $(D^{(N)}_0)^{-1}
D^{(N)}_1=\pi^{(N)}_{\tau_3}$.  Consider the permutation $\tau_5=
\prod_{n\in N}(a_nc_n)$ on $Q_N$.  Notice that $\tau_5=\tau_4\tau_3\tau_4$
and $\tau_2=\tau_5\tau_0$.  By the above $\pi^{(N)}_{\tau_3},
\pi^{(N)}_{\tau_4}\in G(D^{(N)})$, hence $\pi^{(N)}_{\tau_5}\in
G(D^{(N)})$.  Then $\pi^{(N)}_{\tau_2}E^{(N)}_0=\pi^{(N)}_{\tau_5}
\pi^{(N)}_{\tau_0}E^{(N)}_0=\pi^{(N)}_{\tau_5}D^{(N)}_0\in G(D^{(N)})$.
Since $\tau_2(a_n)=a_n$ and $\tau_2(b_n)=b_n$ for all $n\in N$, it easily
follows that transformations $\pi^{(N)}_{\tau_2}$ and $E^{(N)}_0$ commute.
As $\tau_2$ is the product of commuting permutations of odd orders $2n-1$,
$n\in N$, while $E^{(N)}_0$ is an involution, we have that
$(\pi^{(N)}_{\tau_2}E^{(N)}_0)^m=E^{(N)}_0$, where $m=\prod_{n\in N}
(2n-1)$.  In particular, $E^{(N)}_0$ and $\pi^{(N)}_{\tau_2}$ are contained
in $G(D^{(N)})$.  Now Lemma \ref{union1} implies that $\pi^{(N)}_{\tau_1},
E^{(N)}_1\in G(D^{(N)})$.
\end{proof}

Every word $\xi\in(Q_N^\pm)^*$ is assigned a pattern $v$ (i.e., a word in
the alphabet $\{*,*^{-1}\}$) that is obtained from $\xi$ by substituting
$*$ for each occurrence of letters $q\in Q_N$ and substituting $*^{-1}$ for
each occurrence of letters $q^{-1}$, $q\in Q_N$.  We say that $\xi$ follows
the pattern $v$.

Now we introduce an alphabet $P_N^\pm$ that consists of symbols $*_n$ and
$*_n^{-1}$ for all $n\in N$.  A word over the alphabet $P_N^\pm$ is called
a {\em marked pattern}.  Every word $\xi\in(Q_N^\pm)^*$ is assigned a
marked pattern $v\in(P_N^\pm)^*$ that is obtained from $\xi$ as follows.
For any $n\in N$ we substitute $*_n$ for each occurrence of letters $q\in
Q_n$ in $\xi$ and substitute $*_n^{-1}$ for each occurrence of letters
$q^{-1}$, $q\in Q_n$.  We say that $\xi$ follows the marked pattern $v$.
Clearly, the pattern of $\xi$ is uniquely determined by its marked pattern.
Notice that each letter of the alphabet $P_N^\pm$ corresponds to a
connected component of the Moore diagram of the automaton $U^{(N)}$.  Since
$D^{(N)}$ is the dual automaton of $U^{(N)}$, it easily follows that the
$G(D^{(N)})$ action on $(Q_N^\pm)^*$ preserves marked patterns of words.

A word $\xi=q_1q_2\dots q_k\in(Q_N^\pm)^*$ is called freely irreducible if
none of its two-letter subwords $q_1q_2,q_2q_3,\dots,q_{k-1}q_k$ is of the
form $qq^{-1}$ or $q^{-1}q$, where $q\in Q_N$.  Otherwise $\xi$ is called
freely reducible.

\begin{lemma}\label{union3}
For any nonempty word $v\in(P_N^\pm)^*$ there exists a freely irreducible
word $\xi\in(Q_N^\pm)^*$ such that $v$ is the marked pattern of $\xi$ and
the transformation $U^{(N)}_\xi$ acts nontrivially on the first level of
the rooted binary tree $X^*$.
\end{lemma}

\begin{proof}
For any $n\in N$ let us substitute $a_n$ for each occurrence of $*_n$ in
$v$ and $b_n^{-1}$ for each occurrence of $*_n^{-1}$.  We get a nonempty
word $\xi\in(Q_N^\pm)^*$ that follows the marked pattern $v$.  Now let us
modify $\xi$ by changing its last letter.  If this letter is $a_n$ ($n\in
N$), we change it to $c_n$.  If the last letter of $\xi$ is $b_n^{-1}$, we
change it to $c_n^{-1}$.  This yields another word $\eta\in(Q_N^\pm)^*$
that follows the marked pattern $v$.  By construction, $\xi$ and $\eta$ are
freely irreducible.  Furthermore, $U^{(N)}_\eta=A^{(n)}_{c_n}
(A^{(n)}_{a_n})^{-1}U^{(N)}_\xi$ if the last letter of $v$ is $*_n$, $n\in
N$ while $U^{(N)}_\eta=(A^{(n)}_{c_n})^{-1}A^{(n)}_{b_n}U^{(N)}_\xi$ if the
last letter of $v$ is $*_n^{-1}$.  For any $n\in N$ both $A^{(n)}_{c_n}
(A^{(n)}_{a_n})^{-1}$ and $(A^{(n)}_{c_n})^{-1}A^{(n)}_{b_n}$ interchange
one-letter words $0$ and $1$.  It follows that one of the transformations
$U^{(N)}_\xi$ and $U^{(N)}_\eta$ also acts nontrivially on the first level
of the rooted tree $\{0,1\}^*$.
\end{proof}

Given a nonempty, freely irreducible word $\xi\in(Q_N^\pm)^*$, let
$Z_N(\xi)$ denote the set of all freely irreducible words in $(Q_N^\pm)^*$
that follow the same marked pattern as $\xi$ and match $\xi$ completely or
except for the last letter.  Obviously, $\xi\in Z_N(\xi)$, and $\eta\in
Z_N(\xi)$ if and only if $\xi\in Z_N(\eta)$.

\begin{lemma}\label{union4}
For any nonempty word $v\in(P_N^\pm)^*$ there exists a freely irreducible
word $\xi\in(Q_N^\pm)^*$ such that $v$ is the marked pattern of $\xi$ and
the set $Z_N(\xi)$ is contained in one orbit of the $G(D^{(N)})$ action on
$(Q_N^\pm)^*$.
\end{lemma}

\begin{proof}
Let $\tQ_N^\pm=\bigcup_{n\in N}\{a_n,b_n,c_n,a_n^{-1},b_n^{-1},c_n^{-1}\}$.
The set $(\tQ_N^\pm)^*$ of words in the alphabet $\tQ_N^\pm$ is a submonoid
of $(Q_N^\pm)^*$.  Let $h_N:(\tQ_N^\pm)^*\to(Q^\pm)^*$ be the homomorphism
of monoids such that $h_N(a_n)=a$, $h_N(b_n)=b$, $h_N(c_n)=c$,
$h_N(a_n^{-1})=a^{-1}$, $h_N(b_n^{-1})=b^{-1}$, $h_N(c_n^{-1})=c^{-1}$ for
all $n\in N$.  For any $\zeta\in(\tQ_N^\pm)^*$ the word $h_N(\zeta)$
follows the same pattern as $\zeta$.  The word $\zeta$ is uniquely
determined by $h_N(\zeta)$ and the marked pattern of $\zeta$.  If
$h_N(\zeta)$ is freely irreducible then so is $\zeta$ (however $h_N(\zeta)$
can be freely reducible even if $\zeta$ is freely irreducible).  It is easy
to see that $E_0(h_N(\zeta))=h_N(E^{(N)}_0(\zeta))$,
$\pi_{(ab)}(h_N(\zeta))=h_N(\pi^{(N)}_{\si_1}(\zeta))$, and
$\pi_{(bc)}(h_N(\zeta))=h_N(\pi^{(N)}_{\si_2}(\zeta))$, where $\si_1=
\prod_{n\in N}(a_nb_n)$ and $\si_2=\prod_{n\in N}(b_nc_n)$ are permutations
on $Q_N$.  By Proposition \ref{a2}, the group $G(D)$ is generated by $E_0$,
$\pi_{(ab)}$, and $\pi_{(bc)}$.  On the other hand, $E^{(N)}_0,
\pi^{(N)}_{\si_1},\pi^{(N)}_{\si_2}\in G(D^{(N)})$ due to Proposition
\ref{union2}.  Let $\tG$ denote the subgroup of $G(D^{(N)})$ generated by
$E^{(N)}_0$, $\pi^{(N)}_{\si_1}$, and $\pi^{(N)}_{\si_2}$.  It follows that
for any $g_0\in G(D)$ there exists $g\in\tG$ such that $g_0(h_N(\zeta))=
h_N(g(\zeta))$ for all $\zeta\in(\tQ_N^\pm)^*$.  Now Proposition \ref{a5}
implies that words $\zeta_1,\zeta_2\in(\tQ_N^\pm)^*$ are in the same orbit
of the $G(D^{(N)})$ action on $(Q_N^\pm)^*$ whenever they follow the same
marked pattern and the words $h_N(\zeta_1)$, $h_N(\zeta_2)$ are freely
irreducible.

Given a nonempty marked pattern $v\in(P_N^\pm)^*$, let $v_0$ be the word
obtained by deleting the last letter of $v$.  For any $n\in N$ we
substitute $a_n$ for each occurrence of $*_n$ in $v_0$ and $b_n^{-1}$ for
each occurrence of $*_n^{-1}$.  This yields a word $\eta\in(Q_N^\pm)^*$
that follows the marked pattern $v_0$.  Now let $\xi=\eta c_n$ if the last
letter of $v$ is $*_n$, $n\in N$ and let $\xi=\eta c_n^{-1}$ if the last
letter of $v$ is $*_n^{-1}$.  Clearly, $\xi$ is a freely irreducible word
following the marked pattern $v$.  Moreover, $\xi\in(\tQ_N^\pm)^*$ and the
word $h_N(\xi)$ is also freely irreducible.

We shall show that the set $Z_N(\xi)$ is contained in the orbit of $\xi$
under the $G(D^{(N)})$ action on $(Q_N^\pm)^*$.  Take any $\zeta\in
Z_N(\xi)$.  If $\zeta$ is a word over the alphabet $\tQ_N^\pm$ and
$h_N(\zeta)$ is freely irreducible, then it follows from the above that
$\zeta=g(\xi)$ for some $g\in\tG\subset G(D^{(N)})$.  On the other hand,
suppose that the last letter of $\zeta$ is $q_{ni}$ or $q_{ni}^{-1}$, where
$n\in N$, $1\le i\le 2n-2$.  In this case we have $\zeta=
(\pi^{(N)}_\tau)^i(\xi)$, where $\tau=\prod_{n\in N}(c_nq_{n1}\dots
q_{n,2n-2})$.  By Proposition \ref{union2}, $\pi^{(N)}_\tau\in G(D^{(n)})$.

It remains to consider the case when the last letter of $\zeta$ belongs to
$\tQ_N^\pm$ but the word $h_N(\zeta)$ is freely reducible.  There is at
most one $\zeta\in Z_N(\xi)$ with such properties.  It exists if the last
two letters of $v$ are of the form $*_l*_m^{-1}$ or $*_l^{-1}*_m$, where
$l,m\in N$, $l\ne m$.  Assume this is the case.  Then the last letter of
the word $\eta$ is either $a_l$ or $b_l^{-1}$.  Let us change this letter
to $c_l$ or $c_l^{-1}$, respectively.  The resulting word $\eta_1$ follows
the marked pattern $v_0$.  Also, the words $h_N(\eta)$ and $h_N(\eta_1)$
are freely irreducible.  By Proposition \ref{a5}, $h_N(\eta_1)=
g_1(h_N(\eta))$ for some $g_1\in G(D)$.  There exists a unique $\zeta_1\in
(\tQ_N^\pm)^*$ such that $h_N(\zeta_1)=g_1(h_N(\zeta))$ and $v$ is the
marked pattern of $\zeta_1$.  By the above there exists $\tilde g_1\in\tG$
such that $\tilde g_1(\eta)=\eta_1$ and $\tilde g_1(\zeta)=\zeta_1$.  Since
the word $h_N(\zeta)$ is freely reducible, so is $h_N(\zeta_1)$.  On the
other hand, the word $h_N(\eta_1)$, which can be obtained by deleting the
last letter of $h_N(\zeta_1)$, is freely irreducible.  It follows that the
last two letters of $h_N(\zeta_1)$ are $cc^{-1}$ or $c^{-1}c$.  Then the
last two letters of $\zeta_1$ are $c_lc_m^{-1}$ or $c_l^{-1}c_m$.  If
$2m-1$ does not divide $2l-1$ then the word
$(\pi^{(N)}_\tau)^{2l-1}(\zeta_1)$ matches $\zeta_1$ except for the last
letter.  Consequently, the word $\zeta'=\tilde g_1^{-1}
(\pi^{(N)}_\tau)^{2l-1}\tilde g_1(\zeta)$ matches $\zeta$ except for the
last letter.  Since the $G(D^{(N)})$ action preserves marked patterns, the
word $\zeta'$ follows the marked pattern $v$.  Hence $\zeta'\in Z_N(\xi)$.
As $\zeta'\ne\zeta$, it follows from the above that $\zeta'=g(\xi)$ for
some $g\in G(D^{(N)})$.  Then $\zeta=g_0(\xi)$, where $g_0=\tilde g_1^{-1}
(\pi^{(N)}_\tau)^{1-2l}\tilde g_1g\in G(D^{(N)})$.

Now suppose that $2m-1$ divides $2l-1$.  Then $(\pi^{(N)}_\tau)^{2l-1}
(\zeta_1)=\zeta_1$ and the above argument does not apply.  Recall that the
last two letters of $h_N(\zeta_1)$ are $cc^{-1}$ or $c^{-1}c$.  If these
letters are preceded by $b^{-1}$, we let $\zeta_2=\pi^{(N)}_{\si_1}
(\zeta_1)$.  Otherwise they are preceded by $a$ or $h_N(\zeta_1)$ has
length $2$.  In this case, we let $\zeta_2=\zeta_1$.  Further, consider the
permutation $\tau_1=\tau^{2m-1}\si_2\tau^{-(2m-1)}\si_2\tau^{2m-1}$ on
$Q_N$.  Since $\pi^{(N)}_\tau,\pi^{(N)}_{\si_2}\in G(D^{(n)})$, we have
that $\pi^{(N)}_{\tau_1}=(\pi^{(N)}_\tau)^{2m-1}\pi^{(N)}_{\si_2}
(\pi^{(N)}_\tau)^{1-2m}\pi^{(N)}_{\si_2}(\pi^{(N)}_\tau)^{2m-1}\in
G(D^{(N)})$.  It is easy to see that $\tau_1(c_m)=c_m$ and $\tau_1(a_n)=
a_n$ for all $n\in N$.  Since $2m-1<2l-1$, we have $\tau_1(c_l)=b_l$.
Also, for any $n\in N$ we have $\tau_1(b_n)=b_n$ if $2n-1$ divides $2m-1$
and $\tau_1(b_n)=c_n$ otherwise.  It follows that $\zeta_3=
\pi^{(N)}_{\tau_1}(\zeta_2)$ is a word in the alphabet $\tQ_N^\pm$ such
that $h_N(\zeta_3)$ is freely irreducible.  Since $\zeta_3$ follows the
marked pattern $v$, we obtain that $\zeta_3$ belongs to the orbit of $\xi$
under the $G(D^{(N)})$ action.  So does the word $\zeta$.
\end{proof}

\begin{proposition}\label{union5}
Suppose $\xi\in(Q_N^\pm)^*$ is a freely irreducible word.  Then the orbit
of $\xi$ under the action of the group $G(D^{(N)})$ on $(Q_N^\pm)^*$
consists of all freely irreducible words following the same marked pattern
as $\xi$.
\end{proposition}

\begin{theorem}\label{union6}
The group $G(A^{(N)})$ is the free non-Abelian group on generators
$A^{(N)}_q$, $q\in Q_N$.
\end{theorem}

Proposition \ref{union5} is derived from Lemma \ref{union4} in the same way
as Proposition \ref{series6} was derived from Lemma \ref{series5}.  Then
Theorem \ref{union6} is derived from Proposition \ref{union5} and Lemma
\ref{union3} in the same way as Theorem \ref{series7} was derived from
Proposition \ref{series6} and Lemma \ref{series4}.  We omit both proofs.

\section{The Bellaterra automaton and its series}\label{b}

In this section we consider the Bellaterra automaton, a series of automata
of Bellaterra type, and their disjoint unions.  We use the notation of
Sections \ref{a}, \ref{series}, and \ref{union}.

The Bellaterra automaton $B$ is an automaton over the alphabet $X=\{0,1\}$
with the set of internal states $Q=\{a,b,c\}$.  The state transition
function $\hphi$ and the output function $\hpsi$ of $B$ are defined as
follows: $\hphi(a,0)=\hphi(b,1)=c$, $\hphi(a,1)=\hphi(b,0)=b$, $\hphi(c,0)=
\hphi(c,1)=a$; $\hpsi(a,0)=\hpsi(b,0)=\hpsi(c,1)=0$, $\hpsi(a,1)=
\hpsi(b,1)=\hpsi(c,0)=1$.  The Moore diagram of $B$ is depicted in Figure
\ref{fig2}.  It is easy to verify that the inverse automaton of $B$
coincides with $B$.  Besides, the reverse automaton of $B$ can be obtained
from $B$ by renaming its states $a$ and $c$ to $c$ and $a$, respectively.
Lemma \ref{auto6} implies that $B$ is bi-reversible.

The Bellaterra automaton $B$ is closely related to the Aleshin automaton
$A$.  Namely, the two automata share the alphabet, the set of internal
states, and the state transition function.  On the other hand, the output
function $\hpsi$ of $B$ never coincides with the output function $\psi$ of
$A$, that is, $\hpsi(q,x)\ne\psi(q,x)$ for all $q\in Q$ and $x\in X$.

For any integer $n\ge1$ we define a Bellaterra type automaton $B^{(n)}$ as
the automaton that is related to the Aleshin type automaton $A^{(n)}$ in
the same way as the automaton $B$ is related to $A$.  To be precise,
$B^{(n)}$ is an automaton over the alphabet $X=\{0,1\}$ with the set of
states $Q_n$.  The state transition function of $B^{(n)}$ coincides with
that of $A^{(n)}$.  The output function $\hpsi_n$ of $B^{(n)}$ is defined
so that for any $x\in X$ we have $\hpsi_n(q,x)=x$ if $q\in\{a_n,b_n\}$ and
$\hpsi_n(q,x)=1-x$ if $q\in Q_n\setminus\{a_n,b_n\}$.  Then $\hpsi_n(q,x)=
1-\psi_n(q,x)$ for all $q\in Q_n$ and $x\in X$, where $\psi_n$ is the
output function of $A^{(n)}$.  Note that the automaton $B^{(1)}$ coincides
with $B$ up to renaming of the internal states.

In addition, we define a Bellaterra type automaton $B^{(0)}$.  This is an
automaton over the alphabet $X$ with the set of internal states $Q_0$
consisting of a single element $c_0$.  The state transition function
$\hphi_0$ and the output function $\hpsi_0$ of $B^{(0)}$ are defined as
follows: $\hphi_0(c_0,0)=\hphi_0(c_0,1)=c_0$; $\hpsi_0(c_0,0)=1$,
$\hpsi_0(c_0,1)=0$.

It is easy to see that each Bellaterra type automaton $B^{(n)}$ coincides
with its inverse automaton.  The reverse automaton of $B^{(0)}$ coincides
with $B^{(0)}$ as well.  In the case $n\ge1$, the reverse automaton of
$B^{(n)}$ can be obtained from $B^{(n)}$ by renaming its states
$c_n,q_{n1},\dots,q_{n,2n-2},a_n$ to $a_n,q_{n,2n-2},\dots,q_{n1},c_n$,
respectively.  Lemma \ref{auto6} implies that each $B^{(n)}$ is
bi-reversible.

\begin{figure}[t]
%
\includegraphics{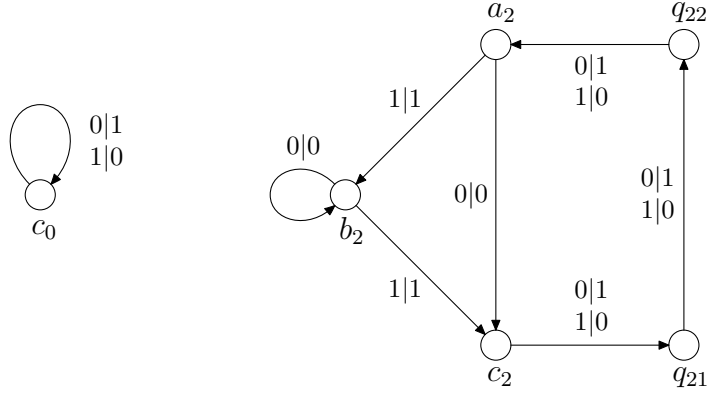}
\includegraphics{second.eps}
\begin{picture}(0,0)(0,85)
\put(35,0)
{
 \begin{picture}(0,0)(0,0)
 \put(33,-15){$c_0$}
 \put(50,17){\footnotesize $\begin{array}{c} 0|1 \\ 1|0 \end{array}$}
 \end{picture}
}
\put(150.5,0)
{
 \begin{picture}(0,0)(0,0)
 \put(90,67){$a_2$}
 \put(34,-17){$b_2$}
 \put(90,-71){$c_2$}
 \put(160,-71){$q_{21}$}
 \put(159.5,67.5){$q_{22}$}
 \put(52.5,33){\footnotesize $1|1$}
 \put(52.5,-39){\footnotesize $1|1$}
 \put(14,15){\footnotesize $0|0$}
 \put(77.5,-3){\footnotesize $0|0$}
 \put(143,-3){\footnotesize $\begin{array}{c} 0|1 \\ 1|0 \end{array}$}
 \put(118,-45.5){\footnotesize $\begin{array}{c} 0|1 \\ 1|0 \end{array}$}
 \put(118,39.5){\footnotesize $\begin{array}{c} 0|1 \\ 1|0 \end{array}$}
 \end{picture}
}
\end{picture}
\vspace{162bp}
\caption{
\label{fig7}
Automaton $B^{(\{0,2\})}$.
}
\end{figure}

Let $N$ be a nonempty set of nonnegative integers.  We denote by $B^{(N)}$
the disjoint union of automata $B^{(n)}$, $n\in N$.  Then $B^{(N)}$ is an
automaton over the alphabet $X=\{0,1\}$ with the set of internal states
$Q_N=\bigcup_{n\in N}Q_n$.  It is bi-reversible since each $B^{(n)}$ is
bi-reversible.  If $0\notin N$, then the automaton $B^{(N)}$ shares its
alphabet, its internal states, and its state transition function with the
automaton $A^{(N)}$ while the output functions of these automata never
coincide.

The relation between automata of Aleshin type and of Bellaterra type
induces a relation between transformations defined by automata of these two
types.

\begin{lemma}\label{b1}
Let $h=B^{(0)}_{c_0}$.  Then

(i) $A_q=hB_q$ and $B_q=hA_q$ for any $q\in\{a,b,c\}$;

(ii) $A^{(n)}_q=hB^{(n)}_q$ and $B^{(n)}_q=hA^{(n)}_q$ for any $n\ge1$ and
$q\in Q_n$;

(iii) $A^{(N)}_q=hB^{(N)}_q$ and $B^{(N)}_q=hA^{(N)}_q$ for any nonempty
set $N$ of positive integers and any $q\in Q_N$.
\end{lemma}

\begin{proof}
The transformation $h$ is the automorphism of the free monoid $\{0,1\}^*$
that interchanges the free generators $0$ and $1$.  For any $w\in X^*$ the
word $h(w)$ can be obtained from $w$ by changing all letters $0$ to $1$ and
all letters $1$ to $0$.

Suppose $\widetilde A$ and $\widetilde B$ are two automata over the
alphabet $X$ such that their sets of internal states and state transition
functions are the same but their output functions never coincide.  It is
easy to see that $\widetilde A_q=h\widetilde B_q$ and $\widetilde B_q=
h\widetilde A_q$ for any internal state $q$ of the automata $\widetilde A$
and $\widetilde B$.  The lemma follows.
\end{proof}

\begin{proposition}\label{b2}
(i) The group $G(A)$ is an index $2$ subgroup of $G(B^{(\{0,1\})})$;

(ii) for any $n\ge1$ the group $G(A^{(n)})$ is an index $2$ subgroup of
$G(B^{(\{0,n\})})$;

(iii) for any nonempty set $N$ of positive integers the group $G(A^{(N)})$
is an index $2$ subgroup of $G(B^{(N\cup\{0\})})$.
\end{proposition}

\begin{proof}
Note that the statement (i) is a particular case of the statement (ii) as
$G(A)=G(A^{(1)})$.  Furthermore, the statement (ii) is a particular case of
the statement (iii) since $A^{(n)}=A^{(\{n\})}$ for any integer $n\ge1$.

Suppose $N$ is a nonempty set of positive integers.  The group $G(A^{(N)})$
is generated by transformations $A^{(N)}_q$, $q\in Q_N$.  The group
$G(B^{(N\cup\{0\})})$ is generated by transformations $h=B^{(0)}_{c_0}$ and
$B^{(N)}_q$, $q\in Q_N$.  By Lemma \ref{b1}, $A^{(N)}_q=hB^{(N)}_q$ and
$B^{(N)}_q=hA^{(N)}_q$ for any $q\in Q_N$.  It follows that the group
$G(B^{(N\cup\{0\})})$ is generated by transformations $h$ and $A^{(N)}_q$,
$q\in Q_N$.  In particular, $G(A^{(N)})\subset G(B^{(N\cup\{0\})})$.

For any $n\ge0$ the automaton $B^{(n)}$ coincides with its inverse.  Lemma
\ref{auto2} implies that $h^2=1$ and $(B^{(N)}_q)^2=1$, $q\in Q_N$.  Then
$hA^{(N)}_qh^{-1}=B^{(N)}_qh=(A^{(N)}_q)^{-1}$ for any $q\in Q_N$.  It
follows that $G(A^{(N)})$ is a normal subgroup of $G(B^{(N\cup\{0\})})$.
Since $h^2=1$, the index of the group $G(A^{(N)})$ in $G(B^{(N\cup\{0\})})$
is at most $2$.  On the other hand, $G(A^{(N)})\ne G(B^{(N\cup\{0\})})$ as
$G(B^{(N\cup\{0\})})$ contains a nontrivial involution $h$ while
$G(A^{(N)})$ is a free group due to Theorem \ref{union6}.  Thus
$G(A^{(N)})$ is an index $2$ subgroup of $G(B^{(N\cup\{0\})})$.
\end{proof}

The relation between groups defined by automata of Aleshin type and of
Bellaterra type allows us to establish the structure of the groups defined
by automata of the latter type.  As the following two theorems show, these
groups are free products of groups of order $2$.

\begin{theorem}[\cite{N}]\label{b3}
The group $G(B)$ is freely generated by involutions $B_a$, $B_b$, $B_c$.
\end{theorem}

\begin{theorem}\label{b4}
(i) For any $n\ge1$ the group $G(B^{(n)})$ is freely generated by $2n+1$
involutions $B^{(n)}_q$, $q\in Q_n$;

(ii) for any nonempty set $N$ of nonnegative integers the group
$G(B^{(N)})$ is freely generated by involutions $B^{(N)}_q$, $q\in Q_N$.
\end{theorem}

To prove Theorems \ref{b3} and \ref{b4}, we need the following lemma.

\begin{lemma}\label{b5}
Suppose that a group $G$ is generated by elements $g_0,g_1,\dots,g_k$
($k\ge1$) of order at most $2$.  Let $H$ be the subgroup of $G$ generated
by elements $h_i=g_0g_i$, $1\le i\le k$.  Then $G$ is freely generated by
$k+1$ involutions $g_0,g_1,\dots,g_k$ if and only if $H$ is the free group
on $k$ generators $h_1,\dots,h_k$.
\end{lemma}

\begin{proof}
Consider an element $h=h_{i_1}^{\eps_1}h_{i_2}^{\eps_2}\dots
h_{i_l}^{\eps_l}$, where $l\ge1$, $1\le i_j\le k$, $\eps_j\in\{-1,1\}$, and
$\eps_j=\eps_{j+1}$ whenever $i_j=i_{j+1}$.  Since $h_i=g_0g_i$ and
$h_i^{-1}=g_ig_0$ for $1\le i\le k$, and $g_0^2=1$, we obtain that
$h=g'_0g_{i_1}g'_1\dots g_{i_l}g'_l$, where each $g'_j$ is equal to $g_0$
or $1$.  Moreover, $g'_j=g_0$ whenever $\eps_j=\eps_{j+1}$.  In particular,
$h\ne1$ if $G$ is freely generated by involutions $g_0,g_1,\dots,g_k$.  It
follows that $H$ is the free group on generators $h_1,\dots,h_k$ if $G$ is
freely generated by involutions $g_0,g_1,\dots,g_k$.

Now assume that $H$ is the free group on generators $h_1,\dots,h_k$.  Then
each $h_i$ has infinite order.  Since $h_i=g_0g_i$ and $g_0^2=g_i^2=1$, it
follows that $g_0\ne1$ and $g_i\ne1$.  Hence each of the elements
$g_0,g_1,\dots,g_k$ has order $2$.  In particular, none of these elements
belongs to the free group $H$.

The group $G$ is freely generated by involutions $g_0,g_1,\dots,g_k$ if
$g\ne1$ for any $g=g_{i_1}\dots g_{i_l}$ such that $l\ge1$, $0\le i_j\le
k$, and $i_j\ne i_{j+1}$.  First consider the case when $l$ is even.  Note
that $g_ig_j=h_i^{-1}h_j$ for $0\le i,j\le n$, where by definition $h_0=1$.
Therefore $g=h_{i_1}^{-1}h_{i_2}\dots h_{i_{l-1}}^{-1}h_{i_l}\in H$.  Since
$h_0=1$, the sequence $h_{i_1}^{-1},h_{i_2},\dots,h_{i_{l-1}}^{-1},h_{i_l}$
can contain the unit elements.  After removing all of them, we obtain a
nonempty sequence in which neighboring elements are not inverses of each
other.  Since $h_1,\dots,h_k$ are free generators, we conclude that
$g\ne1$.  In the case when $l$ is odd, it follows from the above that
$g=g_{i_1}h$, where $h\in H$.  Since $g_{i_1}\notin H$, we have that
$g\notin H$, in particular, $g\ne1$.
\end{proof}

\begin{proofof}{Theorems \ref{b3} and \ref{b4}}
First we observe that Theorem \ref{b3} is a particular case of Theorem
\ref{b4} since the automata $B$ and $B^{(1)}$ coincide up to renaming of
their internal states.  Further, the statement (i) of Theorem \ref{b4} is a
particular case of the statement (ii) since $B^{(n)}=B^{(\{n\})}$ for any
$n\ge1$.

Suppose $N$ is a nonempty set of nonnegative integers such that $0\in N$.
For any $n\in N$ the automaton $B^{(n)}$ coincides with its inverse.  Lemma
\ref{auto2} implies that $(B^{(N)}_q)^2=1$ for all $q\in Q_N$.  If
$N=\{0\}$ then $Q_N=\{c_0\}$ and $G(B^{(N)})$ is a group of order $2$
generated by the involution $h=B^{(0)}_{c_0}$.  Now assume that $N\ne
\{0\}$.  Then $K=N\setminus\{0\}$ is a nonempty set of positive integers.
The group $G(B^{(N)})$ is generated by transformations $h$ and $B^{(K)}_q$,
$q\in Q_K$.  All generators are of order at most $2$.  The group
$G(A^{(K)})$ is the free group on generators $A^{(K)}_q$, $q\in Q_K$ due to
Theorem \ref{union6}.  By Lemma \ref{b1}, $A^{(K)}_q=hB^{(K)}_q$ for any
$q\in Q_K$.  Then Lemma \ref{b5} implies that $G(B^{(N)})$ is freely
generated by involutions $h$ and $B^{(K)}_q$, $q\in Q_K$.

Now consider the case when $N$ is a nonempty set of positive integers.  By
the above the group $G(B^{(N\cup\{0\})})$ is freely generated by
involutions $h$ and $B^{(N)}_q$, $q\in Q_N$.  Clearly, this implies that
the group $G(B^{(N)})$ is freely generated by involutions $B^{(N)}_q$,
$q\in Q_N$.
\end{proofof}

Now we shall establish a relation between transformation groups defined by
the Aleshin type and the Bellaterra type automata with the same set of
internal states.

Since $G(A)$ is the free group on generators $A_a$, $A_b$, $A_c$, there is
a unique homomorphism $\Delta:G(A)\to G(B)$ such that $\Delta(A_a)=B_a$,
$\Delta(A_b)=B_b$, $\Delta(A_c)=B_c$.  Likewise, for any $n\ge1$ there is a
unique homomorphism $\Delta_n:G(A^{(n)})\to G(B^{(n)})$ such that
$\Delta_n(A^{(n)}_q)=B^{(n)}_q$ for all $q\in Q_n$.  Also, for any nonempty
set $N$ of positive integers there is a unique homomorphism $\Delta_N:
G(A^{(N)})\to G(B^{(N)})$ such that $\Delta_N(A^{(N)}_q)=B^{(N)}_q$ for all
$q\in Q_N$.

\begin{proposition}\label{b6}
(i) $G(A)\cap G(B)=\{g\in G(A)\mid \Delta(g)=g\}$;

(ii) $G(A)\cap G(B)$ is the free group on generators $B_aB_b$ and
$B_aB_c$;

(iii) $G(A)\cap G(B)$ is an index $2$ subgroup of $G(B)$;

(iv) $A_p^{-1}A_q=B_pB_q$ for all $p,q\in\{a,b,c\}$.
\end{proposition}

\begin{proof}
Let $h=B^{(0)}_{c_0}$.  By Lemma \ref{b1}, $A_q=hB_q$ for all $q\in
\{a,b,c\}$.  Since the inverse automaton of $B$ coincides with $B$, Lemma
\ref{auto2} implies that $B_a^2=B_b^2=B_c^2=1$.  Then for any $p,q\in
\{a,b,c\}$ we have $A_p^{-1}A_q=(hB_p)^{-1}hB_q=B_p^{-1}B_q=B_pB_q$.

It is easy to see that $\{g\in G(A)\mid \Delta(g)=g\}$ is a subgroup of
$G(A)\cap G(B)$.  Let $\tG$ be the group generated by transformations
$B_aB_b$ and $B_aB_c$.  By the above $\Delta(A_a^{-1}A_b)=B_a^{-1}B_b=
B_aB_b=A_a^{-1}A_b$ and $\Delta(A_a^{-1}A_c)=B_a^{-1}B_c=B_aB_c=
A_a^{-1}A_c$.  It follows that $\tG$ is a subgroup of $\{g\in G(A)\mid
\Delta(g)=g\}$.

By Theorem \ref{b3}, the group $G(B)$ is freely generated by involutions
$B_a$, $B_b$, $B_c$.  Then Lemma \ref{b5} implies that $\tG$ is the free
group on generators $B_aB_b$ and $B_aB_c$.  Note that $B_aB_q\in\tG$ for
all $q\in Q$.  Then for any $p,q\in Q$ we have $B_pB_q=(B_aB_p)^{-1}
B_aB_q\in\tG$.  It follows that for any $g\in G(B)$ at least one of the
transformations $g$ and $B_ag$ belongs to $\tG$.  Therefore the index of
$\tG$ in $G(B)$ is at most $2$.

Note that $B_a\notin G(A)$ as $B_a$ is a nontrivial involution while $G(A)$
is a free group.  Hence $G(A)\cap G(B)\ne G(B)$.  Now it follows from the
above that $\tG=\{g\in G(A)\mid \Delta(g)=g\}=G(A)\cap G(B)$ and this is an
index $2$ subgroup of $G(B)$.
\end{proof}

\begin{proposition}\label{b7}
Let $n$ be a positive integer.  Then

(i) $G(A^{(n)})\cap G(B^{(n)})=\{g\in G(A^{(n)})\mid \Delta_n(g)=g\}$;

(ii) $G(A^{(n)})\cap G(B^{(n)})$ is the free group on $2n$ generators
$B^{(n)}_{a_n}B^{(n)}_q$, $q\in Q_n\setminus\{a_n\}$;

(iii) $G(A^{(n)})\cap G(B^{(n)})$ is an index $2$ subgroup of $G(B^{(n)})$;

(iv) $(A^{(n)}_p)^{-1}A^{(n)}_q=B^{(n)}_pB^{(n)}_q$ for all $p,q\in Q_n$.
\end{proposition}

\begin{proposition}\label{b8}
Let $N$ be a nonempty set of positive integers.  Then

(i) $G(A^{(N)})\cap G(B^{(N)})=\{g\in G(A^{(N)})\mid \Delta_N(g)=g\}$;

(ii) for any $n\in N$ the group $G(A^{(N)})\cap G(B^{(N)})$ is the free
group on generators $B^{(N)}_{a_n}B^{(N)}_q$, $q\in Q_N\setminus\{a_n\}$;

(iii) $G(A^{(N)})\cap G(B^{(N)})$ is an index $2$ subgroup of $G(B^{(N)})$;

(iv) $(A^{(N)}_p)^{-1}A^{(N)}_q=B^{(N)}_pB^{(N)}_q$ for all $p,q\in Q_N$.
\end{proposition}

The proofs of Propositions \ref{b7} and \ref{b8} are completely analogous
to the proof of Proposition \ref{b6} and we omit them. 

Now let us consider the dual automata of the Bellaterra automaton and
automata of Bellaterra type.

\begin{figure}[t]
%
\includegraphics{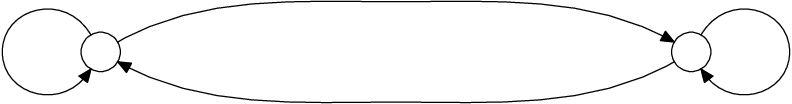}
\begin{picture}(0,0)(-103,45)
\put(-1,-16){$0$}
\put(168,-16){$1$}
\put(79,23){\small $c|a$}
\put(79,-29){\small $c|a$}
\put(-51.7,-3){\small $\begin{array}{r} a|c \\ b|b \end{array}$}
\put(200,-3){\small $\begin{array}{r} a|b \\ b|c \end{array}$}
\end{picture}
\vspace{87bp}
\caption{
\label{fig8}
The dual automaton $\hD$.
}
\end{figure}

Let $\hD$ denote the dual automaton of the Bellaterra automaton $B$.  The
automaton $\hD$ is defined over the alphabet $Q=\{a,b,c\}$, with two
internal states $0$ and $1$.  The Moore diagram of $\hD$ is depicted in
Figure \ref{fig8}.  The automaton $\hD$ is bi-reversible since $B$ is
bi-reversible.

A word $\xi$ over an arbitrary alphabet is called a {\em double letter
word\/} if there are two adjacent letters in $\xi$ that coincide.
Otherwise we call $\xi$ a {\em no-double-letter word}.

The set of no-double-letter words over the alphabet $Q$ forms a subtree of
the rooted ternary tree $Q^*$.  As an unrooted tree, this subtree is
$3$-regular.  However it is not regular as a rooted tree.  The following
proposition shows that the group $G(\hD)$ acts transitively on each level
of the subtree.

\begin{proposition}[\cite{N}]\label{b9}
Suppose $\xi\in Q^*$ is a no-double-letter word.  Then the orbit of $\xi$
under the action of the group $G(\hD)$ on $Q^*$ consists of all
no-double-letter words of the same length as $\xi$.
\end{proposition}

\begin{proof}
Let $\la$ and $\mu$ denote the state transition and output functions of the
automaton $B$.  By $\tilde\la$ and $\tilde\mu$ denote the state transition
and output functions of its dual $\hD$.  Take any $q\in Q$ and $x\in X$.
Let $p=\la(q,x)$ and $y=\mu(q,x)$.  Since $B$ coincides with its inverse
automaton, it follows that $p=\la(q,y)$.  Consequently, $\hD_x(qq)=
\tilde\mu(x,q)\tilde\mu(\tilde\la(x,q),q)=\la(q,x)\la(q,\mu(q,x))=pp$.  It
follows that the set $P=\{qq\mid q\in Q\}\subset Q^*$ is invariant under
$\hD_0$ and $\hD_1$.  Any double letter word $\xi\in Q^*$ is represented as
$\xi_1\xi_0\xi_2$, where $\xi_0\in P$ and $\xi_1,\xi_2\in Q^*$.  For any
$x\in X$ we have $\hD_x(\xi)=\hD_x(\xi_1)\hD_{x_0}(\xi_0)\hD_{x_1}(\xi_2)$,
where $x_0,x_1\in X$.  By the above $\hD_x(\xi)$ is a double letter word.
Thus $\hD_0$ and $\hD_1$ map double letter words to double letter words.
Since these transformations are invertible, they also map no-double-letter
words to no-double-letter words, and so does any $g\in G(\hD)$.

Now we are going to prove that for any no-double-letter words
$\xi_1,\xi_2\in Q^*$ of the same length $l$ there exists $g\in G(\hD)$ such
that $\xi_2=g(\xi_1)$.  The empty word is the only word of length $0$ so it
is no loss to assume that $l>0$.  First consider the case when $l$ is even.
We have $\xi_1=q_1q_2\dots q_{l-1}q_l$ and $\xi_2=p_1p_2\dots p_{l-1}p_l$
for some $q_i,p_i\in Q$, $1\le i\le l$.  Consider two words $\eta_1=
q_1q_2^{-1}\dots q_{l-1}q_l^{-1}$ and $\eta_2=p_1p_2^{-1}\dots
p_{l-1}p_l^{-1}$ over the alphabet $Q^\pm$.  Clearly, $\eta_1$ and $\eta_2$
follow the same pattern.  Furthermore, they are freely irreducible since
$\xi_1$ and $\xi_2$ are no-double-letter words.  By Proposition \ref{a5},
$\eta_2=g_0(\eta_1)$ for some $g_0\in G(D)$.  By Lemma \ref{auto1}, we can
assume that $g_0\in S(D)$. Then $g_0=D_w$ for some word $w\in X^*$.
Proposition \ref{auto4} implies that $U_{\eta_1}(wu)=U_{\eta_1}(w)
U_{\eta_2}(u)$ for any $u\in X^*$.  By Proposition \ref{b6}, $A_p^{-1}A_q=
B_pB_q$ for all $p,q\in Q$.  It follows that $U_{\eta_1}=B_{\xi_1}$ and
$U_{\eta_2}=B_{\xi_2}$.  In particular, $B_{\xi_1}(wu)=B_{\xi_1}(w)
B_{\xi_2}(u)$ for any $u\in X^*$.  Now Proposition \ref{auto4} implies that
$B_{\xi_2}=B_{g(\xi_1)}$, where $g=\hD_w\in G(\hD)$.  By the above
$g(\xi_1)$ is a no-double-letter word.  By Theorem \ref{b3}, the group
$G(B)$ is freely generated by involutions $B_q$, $q\in Q$.  Since $\xi_2$
and $g(\xi_1)$ are no-double-letter words in the alphabet $Q$, the equality
$B_{\xi_2}=B_{g(\xi_1)}$ implies that $\xi_2=g(\xi_1)$.

Now consider the case when $\xi_1$ and $\xi_2$ have odd length.  Obviously,
there exist letters $q_0,p_0\in Q$ such that $\xi_1q_0$ and $\xi_2p_0$ are
no-double-letter words.  Since $\xi_1q_0$ and $\xi_2p_0$ are of the same
even length, it follows from the above that $\xi_2p_0=g(\xi_1q_0)$ for some
$g\in G(\hD)$.  Then $\xi_2=g(\xi_1)$.
\end{proof}

For any integer $n\ge0$ let $\hD^{(n)}$ denote the dual automaton of the
automaton $B^{(n)}$.  The automaton $\hD^{(n)}$ is defined over the
alphabet $Q_n$, with two internal states $0$ and $1$.  It is bi-reversible
since $B^{(n)}$ is bi-reversible.

\begin{proposition}\label{b10}
Let $n\ge1$ and suppose $\xi\in Q_n^*$ is a no-double-letter word.  Then
the orbit of $\xi$ under the action of the group $G(\hD^{(n)})$ on $Q_n^*$
consists of all no-double-letter words of the same length as $\xi$.
\end{proposition}

The proof of Proposition \ref{b10} is completely analogous to the above
proof of Proposition \ref{b9} and we omit it.

\bigskip

{\sc
\begin{raggedright}
Department of Mathematics\\
Texas A\&M University\\
College Station, TX 77843--3368
\end{raggedright}
}


\begin{thebibliography}{BGN}

\bibitem[Ale]{A}
S.\,V.\,Aleshin.  A free group of finite automata. {\em Mosc. Univ. Math.
Bull.} {\bf 38} (1983), no.\,4, 10--13.

\bibitem[BGN]{BGN}
L.\,Bartholdi, R.\,Grigorchuk, V.\,Nekrashevych.  From fractal groups to
fractal sets. {\em Grabner P. (ed.) et al., Fractals in Graz 2001.
Analysis, dynamics, geometry, stochastics. Proceedings of the conference,
Graz, Austria, June 2001}, 25--118.  Trends in Math., Birkh\"auser, Basel,
2003.

\bibitem[BGS]{BGS}
L.\,Bartholdi, R.\,I.\,Grigorchuk, Z.\,\v{S}uni\'k.  Branch groups. {\em
Hazewinkel M. (ed.), Handbook of Algebra, vol.\,3}, 989--1112. Elsevier,
Amsterdam, 2003.

\bibitem[BS]{BS}
A.\,M.\,Brunner, S.\,Sidki.  The generation of
$\mathrm{GL}(n,\bZ)$ by finite state automata. {\em Int. J. Algebra
Comput.} {\bf 8} (1998), no.\,1, 127--139.

\bibitem[Gri]{G}
R.\,I.\,Grigorchuk.  Cancellative semigroups of power growth. {\em Math.
Notes\/} {\bf 43} (1988), no.\,3, 175--183.

\bibitem[GM]{GM}
Y.\,Glasner, S.\,Mozes.  Automata and square complexes. {\em Geom.
Dedicata\/} {\bf 111} (2005), 43--64.

\bibitem[GNS]{GNS}
R.\,I.\,Grigorchuk, V.\,V.\,Nekrashevich, V.\,I.\,Sushchanskii.  Automata,
dynamical systems, and groups. {\em Grigorchuk R.\,I. (ed.), Dynamical
systems, automata, and infinite groups.} {\em Proc. Steklov Inst. Math.}
{\bf 231} (2000), 128--203.

\bibitem[Ho\v{r}]{H}
J.\,Ho\v{r}ej\v{s}.  Transformations defined by finite automata. {\em
Probl. Kibernetiki\/} {\bf 9} (1963), 23--26 (in Russian).

\bibitem[MNS]{MNS}
O.\,Macedo\'nska, V.\,Nekrashevich, V.\,Sushchanskij.  Commensurators of
groups and reversible automata. {\em Dopov. Nats. Akad. Nauk Ukr., Mat.
Pryr. Tekh. Nauky\/} (2000), no.\,12, 36--39.

\bibitem[Nek]{N}
V.\,Nekrashevych.  {\em Self-similar groups}. Math. Surveys and Monographs
{\bf 117}. Amer. Math. Soc., Providence, RI, 2005.

\bibitem[Oli1]{O1}
A.\,S.\,Olijnyk.  Free groups of automaton permutations. {\em Dopov. Nats.
Akad. Nauk Ukr., Mat. Pryr. Tekh. Nauky\/} (1998), no.\,7, 40--44 (in
Ukrainian).

\bibitem[Oli2]{O2}
A.\,S.\,Olijnyk.  Free products of finite groups and groups of finitely
automatic permutations. {\em Grigorchuk R.\,I. (ed.), Dynamical
systems, automata, and infinite groups.} {\em Proc. Steklov Inst. Math.}
{\bf 231} (2000), 308--315.

\bibitem[Sid]{S}
S.\,Sidki.  Automorphisms of one-rooted trees: growth, circuit structure,
and acyclicity. {\em J. Math. Sci., NY\/} {\bf 100} (2000), no.\,1,
1925--1943.

\bibitem[VV]{VV}
M.\,Vorobets, Y.\,Vorobets.  On a free group of transformations defined by
an automaton.  Preprint, 2006 (arXiv:math.GR/0601231).

\end{thebibliography}
\end{document}